\documentclass[10pt]{article}
\pdftrailerid{}
\usepackage[utf8]{inputenc}
\usepackage[T1]{fontenc}
\usepackage{amsmath,amssymb,dsfont}
\numberwithin{equation}{section}
\usepackage{microtype}
\usepackage{graphicx,tikz,pgfplots}
\graphicspath{{images/}}
\pgfplotsset{compat=newest}
\usepackage[hyperref,amsmath,thmmarks]{ntheorem}
\usepackage{aliascnt}
\usepackage[a4paper,centering,bindingoffset=0cm,marginpar=2cm,margin=2.5cm]{geometry}
\usepackage[pagestyles]{titlesec}
\usepackage[font=footnotesize,format=plain,labelfont=sc,textfont=sl,width=0.75\textwidth,labelsep=period]{caption}

\makeatletter
\begingroup
\endlinechar=-1\everyeof{\noexpand}
\edef\x{\endgroup\def\noexpand\bibinputs{\@@input|"kpsewhich infmath.bib" }}\x
\makeatother
\ifx\bibinputs\empty\else\def\havecscbib{}\fi

\usepackage[backend=biber,maxnames=10,backref=true,hyperref=true,giveninits=true,safeinputenc]{biblatex}
\ifdefined\havecscbib
    \ifdefined\extract
    \else
    \fi
\else
    \ifdefined\extract
      \errmessage{cannot extract without detected csc bibliography}
    \fi
\fi

\DefineBibliographyStrings{english}{%
	backrefpage = {cited on page},
	backrefpages = {cited on pages},
}

\DeclareFieldFormat[report]{title}{``#1''}
\DeclareFieldFormat[book]{title}{``#1''}
\AtEveryBibitem{\clearfield{url}}
\AtEveryBibitem{\clearfield{note}}

\titleformat{\section}[block]{\large\sc\filcenter}{\thesection.}{0.5ex}{}[]
\titleformat{\subsection}[runin]{\bf}{\thesubsection.}{0.5ex}{}[.]

\usepackage[pdftex,colorlinks=true,linkcolor=blue,citecolor=green,urlcolor=blue,bookmarks=true,bookmarksnumbered=true]{hyperref}
\hypersetup
{
    bookmarksnumbered,
    breaklinks,
    pdfborderstyle=,
    colorlinks,
    citecolor=[rgb]{0,.65,0},
    linkcolor=[rgb]{0,0,.5},
    urlcolor=[rgb]{0,0,.5},
    pdfauthor={Authors},
    pdfsubject={Subject},
    pdftitle={Title},
    pdfkeywords={Keywords}
}

\newpagestyle{headers}
{
	\headrule
	\sethead[\footnotesize\thepage][\footnotesize\sc Authors][]{}{\footnotesize\sc Computational inverse scattering with internal sources: a RKHS approach}{\footnotesize\thepage}
	\setfoot{}{}{}
}
\newaliascnt{proposition}{lemma}

\aliascntresetthe{proposition}

\newaliascnt{corollary}{lemma}

\aliascntresetthe{corollary}

\newaliascnt{theorem}{lemma}

\aliascntresetthe{theorem}

\theorembodyfont{\normalfont}
\newaliascnt{definition}{lemma}

\aliascntresetthe{definition}

\newaliascnt{assumption}{lemma}

\aliascntresetthe{assumption}

\theoremstyle{nonumberplain}
\theoremseparator{:}
\theoremheaderfont{\normalfont\itshape}

\theoremsymbol{\ensuremath{\square}}

\let\RE\Re
\let\Re=\undefined
\DeclareMathOperator{\Re}{\RE e}
\let\IM\Im
\let\Im=\undefined
\DeclareMathOperator{\Im}{\IM m}

\let\ii\i
\renewcommand{\i}{\mathrm i}

\title{Computational inverse scattering with internal sources: a reproducing kernel Hilbert space approach}
\author{Yakun Dong$^1$\\{\footnotesize\href{mailto:email}{yakun.dong@univie.ac.at}}
	\and Kamran Sadiq$^{2}$\\{\footnotesize\href{mailto:email}{Kamran.Sadiq@oeaw.ac.at}}
	\and Otmar Scherzer$^{1,2,3}$\\{\footnotesize\href{mailto:email}{otmar.scherzer@univie.ac.at}}
	\and John C. Schotland$^{4}$\\{\footnotesize\href{mailto:email}{john.schotland@yale.edu}}}

\date{}
\DeclareMathOperator{\br}{\mathbf{r}}
\DeclareMathOperator*{\argmin}{argmin}
\newcommand{\Laplace}{\mathop{}\!\mathbin\bigtriangleup }

\usepackage{subfigure}

\begin{document}
	
	\maketitle
	\thispagestyle{empty}
	\begin{center}
		\parbox[t]{12em}{\footnotesize
			\hspace*{-1ex}$^1$Faculty of Mathematics\\
			University of Vienna\\
			Oskar-Morgenstern-Platz 1\\
			A-1090 Vienna, Austria}
		\hfil
		\parbox[t]{17em}{\footnotesize
			\hspace*{-1ex}$^2$Johann Radon Institute for Computational\\
			\hspace*{1em}and Applied Mathematics (RICAM)\\
			Altenbergerstraße 69\\
			A-4040 Linz, Austria}
	\end{center}
	
	\begin{center}
		\parbox[t]{19em}{\footnotesize
			\hspace*{-1ex}$^3$Christian Doppler Laboratory for Mathematical  \\
			Modeling and Simulation of Next Generations\\
			of Ultrasound Devices (MaMSi)\\
			Oskar-Morgenstern-Platz 1\\
			A-1090 Vienna, Austria}
		\hfil
		\parbox[t]{17em}{\footnotesize
			\hspace*{-1ex}$^4$Department of Mathematics and \\
			Department of Physics\\
			Yale University\\
			New Haven, CT 06520-8283, USA}
	\end{center}

	\newcommand{\red}[1]{{\color{red}{#1}}}
	\newcommand{\commentK}[1]{\textcolor{blue}{\textbf{Kamran:} #1}}
	\newcommand{\commentJ}[1]{\textcolor{green}{\textbf{John:} #1}}
	\newcommand{\commentY}[1]{\textcolor{cyan}{\textbf{Yakun:} #1}}

	
	\section*{Abstract}
	We present a method to reconstruct the dielectric susceptibility (scattering potential) of an inhomogeneous scattering medium, based on the solution to the inverse scattering problem with internal sources.
	We employ the theory of reproducing kernel Hilbert spaces, together with regularization to recover the susceptibility of two- and three-dimensional scattering media. Numerical examples illustrate the effectiveness of the proposed reconstruction method.
	
	\section{Introduction} \label{sec:intro}
	
	The inverse scattering problem is concerned with determining the structure of a scattering medium from measurements of the scattered field.  Such problems have numerous applications including medical imaging, nondestructive testing, remote and radar sensing, ocean acoustics and geophysical exploration; see e.g. 
	\cite{Kuc13, Eps03, WegAraCarStoMat03, CarScho21}, and the references therein.
	
	In this paper we investigate an inverse scattering problem with internal sources, and show that it is possible to recover the dielectric susceptibility (scattering potenial) of an inhomogeneous scattering medium.
	Suitable internal sources take the form of photoactivated fluorophores, which are used in superresolution microscopy \cite{BetPatSouLinOle06,HelSahBatZhuHei15}, and 
	subwavelength bubbles, which are used in ultrasound localization microscopy \cite{ErrPiePezDesLen15}.
	The considered inverse problem is ill-posed. We apply the theory of reproducing kernel Hilbert spaces (RKHS)  \cite{Aro50b,KimWah70} accompanied by tools from regularization theory (see e.g. \cite{TikArs77, EngHanNeu96})
	to reconstruct the susceptibility.

	In \cite{GilLevScho18}, a similar inverse problem is considered, also using an integral equation approach. The integral equation is converted into a linear system of algebraic equations, which is solved by means of a regularized pseudoinverse to recover the susceptibility.  
	Instead, here we use the RKHS approach, which yields reconstructed images of higher quality compared to \cite{GilLevScho18}.

	This paper is organized as follows. In Section \ref{sec:method}, we formulate the forward problem and the associated inverse scattering problem.  We also briefly recall some tools from the theory of regularization and RKHS, which we need for our reconstruction method.  In addition, we present the method for both two- and three-dimensional inverse problems.
	In Section \ref{sec:numerics}, we present the results obtained by applying the RKHS reconstruction method to three numerical experiments. We also present some comparison results to \cite{GilLevScho18}, as shown in Table \ref{tb:comp-john}.
	Finally, in Section \ref{sec:conclusion}, we present our concluding remarks.

	\section{Inverse scattering problem with internal sources} \label{sec:method}

	We consider an experiment in which light from a monochromatic source propagates in an inhomogeneous medium in two dimensions. The three-dimensional case is taken up in Section \ref{subsec:3d}.
	The source is taken to be located in the interior of the medium and to consist of a single photoactivated fluorescent molecule. The scattered light is then registered by a detector in the far field for each source. 
	We ignore the effects of polarization and consider a scalar model of the optical field, which satisfies the time-harmonic 
	wave equation,
	\begin{align}\label{eq:wave}
		\Laplace U(\br) +k^2 (1+\eta (\br)) U(\br) = - a(\br_1) \delta(\br - \br_1), 
	\end{align}	
	where $U$ is the total field, $\eta(\br)$ is the susceptibility of the medium at the position $\br$, and $k$ is the wave number. Since the source consists of a single molecule, it is
	taken to be a point source with position $\br_1$ and amplitude $a(\br_1)$. The total field may be decomposed into the incident field $U_i$ and scattered field $U_s$. Following standard procedure \cite{CarScho21,BorWol99}, the total field satisfies the integral equation
	\begin{equation}\label{eq:Int_s}
		U(\br ) = U_i(\br ) +k^2 \int G(\br ,\br') \eta(\br') U(\br') d\br ',
	\end{equation} 
	where\begin{equation}
		{G}(\br ,\br') := \frac{\mathrm{i}}{4} H_0^{(1)}(k |\br -\br'|)
	\end{equation}
	is the two-dimensional Green's function, with $H_0^{(1)}$  being the Hankel function of first kind \cite{ColKre98}.
	The Green's function satisfies the equation
	\begin{equation}
		\displaystyle \Laplace_{\mathbf{r}} G(\mathbf{r},\mathbf{r'}) +k^2 G(\mathbf{r},\mathbf{r'}) = - \delta(\mathbf{r} - \mathbf{r'}),
	\end{equation}
	together with the Sommerfeld radiation condition \cite{ColKre98}.
	
	In the far-zone, the Hankel function has the asymptotic form, see e.g. \cite[Section 3.4 (3.59)]{ColKre98}:
	\begin{align}
		\frac{\mathrm{i}}{4} H_0^{(1)}(k R) & \sim \frac{e^{\mathrm{i} \frac{\pi}{2}}}{4} \sqrt{\frac{2}{\pi k R}} e^{-\frac{\pi}{4} \mathrm{i}} e^{\mathrm{i} k R}=\frac{e^{\mathrm{i} \frac{\pi}{4}}}{\sqrt{8 \pi k}} \frac{e^{\mathrm{i} k R}}{\sqrt{R}},
	\end{align}
	which propagates as a cylindrical wave with amplitude $R^{-1 / 2}$. 
	Thus, in the far field,  we obtain the asymptotic formula for the two-dimensional Green's function,
	\begin{align}\nonumber
		G\left(\mathbf{r}, \mathbf{r}^{\prime}\right)&=
		\frac{\mathrm{i}}{4} H_0^{(1)}\left(k\left|\mathbf{r}-\mathbf{r}^{\prime}\right|\right) 
		\sim \frac{e^{\mathrm{i} \frac{\pi}{4}}}{\sqrt{8 \pi k}} \frac{e^{\mathrm{i} k\left|\mathbf{r}-\mathbf{r}^{\prime}\right|}}{\sqrt{\left|\mathbf{r}-\mathbf{r}^{\prime}\right|}} \\ \label{eq:asy-g}
		&\simeq \frac{e^{\mathrm{i} \frac{\pi}{4}}}{\sqrt{8 \pi k}} \frac{e^{\mathrm{i} k|\mathbf{r}|}}{\sqrt{|\mathbf{r}|}} e^{-\mathrm{i} k \hat{\mathbf{r}} \cdot \mathbf{r}^{\prime}}\left(1+\frac{1}{2} \frac{\widehat{\mathbf{r}} \cdot \mathbf{r}^{\prime}}{|\mathbf{r}|}\right).
	\end{align}
	Since the incident field is given by $U_i(\mathbf{r})=a\left(\mathbf{r}_1\right) G\left(\mathbf{r}, \mathbf{r}_1\right)$,  we derive  by using \eqref{eq:asy-g} the asymptotic behavior of the total field,
	\begin{equation}
		\begin{aligned}
			U(\mathbf{r})  &=a\left(\mathbf{r}_1\right) G\left(\mathbf{r}, \mathbf{r}_1\right)+k^2 \int G\left(\mathbf{r}, \mathbf{r}^{\prime}\right) \eta\left(\mathbf{r}^{\prime}\right) U\left(\mathbf{r}^{\prime}\right) d \mathbf{r}^{\prime} \\
			& \sim \frac{e^{\mathrm{i} \frac{\pi}{4}}}{\sqrt{8 \pi}} \frac{e^{\mathrm{i} k|\mathbf{r}|}}{\sqrt{|\mathbf{r}|}} 
			\left(
			e^{-\mathrm{i} k \hat{\mathbf{r}} \cdot \mathbf{r}_1} \frac{a\left(\mathbf{r}_1\right)}{\sqrt{k}}
			+
			k^{\frac{3}{2}} \int e^{-\mathrm{i} k \hat{\mathbf{r}} \cdot \mathbf{r}^{\prime}} \eta\left(\mathbf{r}^{\prime}\right) U\left(\mathbf{r}^{\prime}\right) d \mathbf{r}^{\prime}
			\right)
		\end{aligned}
	\end{equation}
	Thus, in the far zone of the scatterer, the field behaves as an outgoing cylindrical wave:
	\begin{equation}\label{eq:asy-u}
		U(\mathbf{r}) \sim \frac{e^{\mathrm{i} \frac{\pi}{4}}}{\sqrt{8 \pi}} \frac{e^{\mathrm{i} k|\mathbf{r}|}}{\sqrt{|\mathbf{r}|}}\left(\frac{a(\br_1)}{\sqrt{k}} e^{-\mathrm{i} k \hat{\mathbf{r} \cdot \mathbf{r}_1}}+A(\mathbf{r})\right),
	\end{equation}
	where the scattering amplitude $A$ is defined by
	\begin{equation}\label{eq:def-ScatAmp}
		A(\mathbf{r}):=k^{\frac{3}{2}} \int_{\Omega} e^{-\i k \hat{\mathbf{r}} \cdot \mathbf{r}^{\prime}} \eta\left(\mathbf{r}^{\prime}\right) U\left(\mathbf{r}^{\prime}\right) d \mathbf{r}^{\prime},
	\end{equation} and $\Omega$ is the volume of the scatterer.
	
	In the weak-scattering approximation, the scattering amplitude can be calculated by replacing $U$ with $U_i$. Since the scattering amplitude $A$ is recorded in the far zone at at the point $\br_2$, the scattering amplitude can be rewritten as:
	\begin{equation}\label{eq:ScatAmp-r1r2}
		\begin{aligned}
			A(\br_1,\br_2) &= a(\br_1) k^{\frac{3}{2}}  \int_{\Omega} e^{-\i k \widehat{\br}_2 \cdot \br} \eta(\br) G(\br , \br_1 ) d\br,
		\end{aligned}
	\end{equation}
	where the dependence of the scattering amplitude $A(\br_1,\br_2)$ on the source position $\br_1$ and the detector position $\br_2$ has been made explicit. The inverse problem thus consists of solving the integral equation \eqref{eq:ScatAmp-r1r2} to recover the susceptibility $\eta$.
	
	If we differentiate \eqref{eq:ScatAmp-r1r2} with respect to $\br_1$, use the fact that $\displaystyle \Laplace_{\mathbf{r}} G(\mathbf{r},\mathbf{r'}) +k^2 G(\mathbf{r},\mathbf{r'}) = - \delta(\mathbf{r} - \mathbf{r'}) $ for the Green's function, and assume unit amplitude for all sources in \eqref{eq:ScatAmp-r1r2}, we obtain 
	\begin{equation}
		\begin{aligned}\label{eq:diffA}
			\Laplace_{\br_1} A(\br_1,\br_2) + k^2 A(\br_1,\br_2) &=k^{\frac{3}{2}}  \int_{\Omega} e^{-\mathrm{i} k \widehat{\br}_2 \cdot \br} \eta(\br) 
			\left( \Laplace_{\br_1}G(\br , \br_1 )+ k^2 G(\br , \br_1 )\right)  d\br \\
			&=  k^{\frac{3}{2}}  \int_{\Omega} e^{-\mathrm{i} k \widehat{\br}_2 \cdot \br} \eta(\br)  \left(- \delta(\br - \br_1) \right)  d\br  
			=  - k^{\frac{3}{2}}  e^{-\mathrm{i} k \widehat{\br}_2 \cdot \br_1} \eta(\br_1).
		\end{aligned}
	\end{equation}
	We thus obtain an inversion formula for $\eta$ of the form
	\begin{align}\label{eq:eta_2d}
		\eta(\br_1) = -\frac{e^{\mathrm{i} k \widehat{\br}_2 \cdot \br_1}}{k^{\frac{3}{2}}   } \left( \Laplace_{\br_1} A(\br_1,\br_2) + k^2 A(\br_1,\br_2) \right).
	\end{align}
	Theoretically, reconstruction of the susceptibility $\eta$ using \eqref{eq:eta_2d} only requires the position of the sources $\br_1$. However, numerically, reconstructions of $\eta$ based on different observations of the scattering amplitude from different detectors $\br_2$ are obtained. Thus we can determine the susceptibility $\eta(\br_1)$ by averaging the reconstructions:
	\begin{equation}\label{ eq:eq_2d_ave}
		\eta(\br_1) =\frac{\sum_{\br_2}\eta(\br_1;\br_2)}{n_d}, 
	\end{equation}
	where the dependence of $\eta$ on the detector positions is made explicit and $n_d$ represents the number of detectors $\br_2$.
	
	We illustrate the proposed reconstruction method in both two-dimensional and three dimensional scattering media (see subsection \ref{subsec:3d} below).  To avoid any ambiguity, we call the inversion of the susceptibility in a two-dimensional medium as the two-dimensional reconstruction method, while the inversion in a three-dimensional medium is called the three-dimensional reconstruction method. 
	
	

	\subsection{Numerical analysis of the forward problem} We begin with the numerical analysis of the forward problem. Since the scattering amplitude is computed from \eqref{eq:ScatAmp-r1r2} and there is a logarithmic singularity in the two-dimensional Green's function, a careful evaluation of the integral is required. To proceed, we discretize $\Omega$ into pixels $C_k\left(k=1,2, \cdots N_\nu\right) $ of area  $h^2$, where $N_\nu $ is the total number of pixels. Assuming unit amplitude $a=1$ for
	all sources, we split the integral of \eqref{eq:ScatAmp-r1r2} into singular- and non-singular parts:  
	\begin{align}\label{A_int2D}
		k^{\frac{3}{2}} \int_{\Omega} e^{-\i k \widehat{\br}_2 \cdot \br} \eta(\br) G(\br , \br_1 ) d\br=k^{\frac{3}{2}} \sum_{k=1}^{N_\nu} \int_{\Omega_k \backslash \Omega_s} e^{-\i k \widehat{\br}_2 \cdot \br} \eta(\br) G(\br , \br_1 ) d\br+k^{\frac{3}{2}} \int_{\Omega_s} e^{-\i k \widehat{\br}_2 \cdot \br} \eta(\br) G(\br , \br_1 ) d\br,
	\end{align}
	where $\Omega_s$ is the pixel containing the singularity.
	We assume that $k h \ll 1$. As the non-singular part can be computed directly, the singular part remains to be estimated.
	
	The Green's function $G\left(\mathbf{r}, \mathbf{r}^{\prime}\right)$ is weakly singular at $\mathbf{r} = \mathbf{r}^{\prime}$.  
	We will approximate the two-dimensional Green's function near the singularity using its asymptotic form:
	\begin{align}\label{eq:singular-2d-g}
		\begin{aligned}
			G\left(\mathbf{r}, \mathbf{r}^{\prime}\right)=&\frac{1}{2 \pi} \ln \frac{1}{\left|\mathbf{r}-\mathbf{r}^{\prime}\right|}+\frac{\mathrm{i}}{4}
			-\frac{\gamma}{2 \pi}-\frac{1}{2 \pi} \ln \frac{k}{2}+O\left(\left|\mathbf{r}-\mathbf{r}^{\prime}\right|^2 \ln \frac{1}{\left|\mathbf{r}-\mathbf{r}^{\prime}\right|}\right),
		\end{aligned}
	\end{align}
	for $\left|\mathbf{r}-\mathbf{r}^{\prime}\right| \rightarrow 0$, where $\gamma=0.5772156$ is the Euler-Mascheroni constant. Thus, the singular part in \eqref{A_int2D} can be computed as
	\begin{align*}
		\begin{aligned}
			k^{\frac{3}{2}} \int_{\Omega_s} \eta(\mathbf{r}) e^{-\mathrm{i} k \widehat{\mathbf{r}}_2 \cdot \mathbf{r}} G\left(\mathbf{r}, \mathbf{r}_1\right) d \mathbf{r} &\approx
			\frac{k^{\frac{3}{2}}}{2 \pi} \eta(\widetilde{\mathbf{r}}) e^{-\mathrm{i} k \hat{\mathbf{r}}_2 \cdot \tilde{\mathbf{r}}}\left[\left(\frac{\mathrm{i} \pi}{2}-\gamma-\ln \frac{k}{2}\right) h^2+\int_{-\frac{h}{2}}^{\frac{h}{2}} \int_{-\frac{h}{2}}^{\frac{h}{2}}\left(\ln \frac{1}{\left|\mathbf{r}-\mathbf{r}_1\right|}\right) d \mathbf{r}\right] \\
			& =\eta(\widetilde{\mathbf{r}}) e^{-\mathrm{i} k \hat{\mathbf{r}}_2 \cdot \widetilde{\mathbf{r}}} \frac{k^{\frac{3}{2}} h^2}{2 \pi}\left(\xi-\ln \frac{h k}{2}\right),
		\end{aligned}
	\end{align*}
	where  $\displaystyle \xi=\frac{1}{2}(3+\ln 2)-\frac{\pi}{4}-\gamma+\frac{\mathrm{i} \pi}{2}$ and $\widetilde{\mathbf{r}}$ is the central point in $\Omega_s$.
	
	\subsection{Numerical analysis of the inverse problem}
	Next, we focus on estimating the susceptibility through differentiation of the measured scattering amplitude $A(\mathbf{r}_1,\mathbf{r}_2)$ using \eqref{eq:eta_2d}. We apply the theory of RKHS \cite{Aro50b,KimWah70}, to carry out the indicated differentiation, 
	which allows for an explicit formulation of a functional by means of a reproducing kernel \cite{SaiSaw16,Per22}.
	In the setting of our problem, we identify a functional $f$ within a space $\mathcal{W}$ that approximates the scattering amplitude $A(\mathbf{x})$. Here $\mathbf{x} = (\mathbf{r}_1,\mathbf{r}_2) = \mathcal{D} \in \mathbb{R}^4$, where $\mathbf{r}_1$ and $\mathbf{r}_2$ are the position of the sources and detectors, respectively. This leads to the following optimization problem:
	\begin{equation}\label{eq:op-min}
		\argmin_{f} \frac{1}{n}\sum_{i = 1}^{n}|\bar{A}_i-f(\mathbf{x}_i)|^2+ \lambda\|f\|^2_{\mathcal{W}},\quad \bar{A}_i = A(\mathbf{x}_i),
	\end{equation}
	where $n = n_s \times n_d$, with $n_s$ representing the number of sources and $n_d$ representing the number of detectors.
	
	The following result (commonly referred as the representer theorem 
	\cite{KimWah70, KimWah71}) demonstrates that the minimizer of problem \eqref{eq:op-min} still lies in the space $\mathcal{W}$ under certain conditions. 
	This result effectively reduces the original infinite-dimensional minimization problem \eqref{eq:op-min} to a finite-dimensional one for the scalar coefficients. More precisely, if the space $\mathcal{W}$ has a reproducing kernel $\kappa: X\times X\rightarrow \mathbb{R}$, then the minimizer $f^*$ of the problem \eqref{eq:op-min} lies in a finite-dimensional subspace of the infinite-dimensional space $\mathcal{W}$, 
	where
	\begin{equation}\label{eq:f-min}
		f^*(\cdot) = \sum_{i=1}^{n}c^*_i\kappa(\mathbf{x}_i,\cdot), \ \ c^*_i\in\mathbb{R}, \ \ i=1,2,\dots,n,\end{equation}
	and $\mathbf{x}_i$ is the known input.
	
	Let $\mathcal{W}$ be a reproducing kernel Hilbert space,  $\kappa$ be its reproducing kernel, and  $K$ be the Gram matrix, with $ {K}_{i,j}=\kappa(\mathbf{x}_i,\mathbf{x}_j)$, where $\mathbf{x}_i,\mathbf{x}_j \in \mathcal{D}$. Then any functional $f\in \text{span}\{ \kappa(\mathbf{x}_i,\cdot), i = 1,2,\dots,n\}$, can be written as $f=\sum^n_{i=1}c_i\kappa(\mathbf{x}_i,\cdot) $. Thus the vector $(f(\mathbf{x}_1),\dots,f(\mathbf{x}_n)) = \bar{c} {K}$, where $\bar{c} = (c_1,\cdots,c_n)$. The components of \eqref{eq:op-min} can be rewritten in terms of the Gram matrix ${K}$ and the coefficient $\bar{c}$ as
	$$n^{-1}\sum_{i = 1}^{n}|\bar{A}_i-f(\mathbf{x}_i)|^2 = n^{-1}\|\bar{A} - \bar{c}{K} \|^2_{\mathbb{R}^n},$$
	where $\bar{A} = (\bar{A}_1,\dots,\bar{A}_n) $ and 
	\begin{align*}
		\lambda \|f\|^2_{\mathcal{W}} &= \lambda \langle \sum^n_{i=1}c_i\kappa(\mathbf{x}_i,\cdot) , \sum^n_{j=1}c_j\kappa(\mathbf{x}_j,\cdot)\rangle_{\mathcal{W}}
		= \lambda \sum^n_{i,j=1}c_ic_j\kappa(\mathbf{x}_i,\mathbf{x}_j)  = \lambda\langle \bar{c} K,\bar{c}  \rangle_{\mathbb{R}^n},
	\end{align*}
	where we have used the symmetry property of the kernel function. This allows us to consider the optimization problem \eqref{eq:op-min} as a function of $\bar{c}$ admitting a representation,
	\begin{equation}
		\begin{aligned}
			\mathcal{M}(\bar{c})  &:= n^{-1}\sum_{i = 1}^{n}|\bar{A}_i-f(\mathbf{x}_i)|^2+ \lambda\|f\|^2_{\mathcal{W}} =n^{-1}\|\bar{A} - \bar{c}{K} \|^2_{\mathbb{R}^n}+ \lambda \langle \bar{c} K,\bar{c}  \rangle_{\mathbb{R}^n} \\
			&=n^{-1}\|\bar{A}\|^2_{\mathbb{R}^n}  -2 n^{-1}\langle \bar{c} ,\bar{A} K  \rangle_{\mathbb{R}^n} + \langle \bar{c} ,\bar{c}(\lambda K+n^{-1}K^2    ) \rangle_{\mathbb{R}^n},
		\end{aligned}
	\end{equation}
	where in the last equality we use the fact that $K$ is a symmetric matrix. By the representer theorem, if the minimizer of the problem \eqref{eq:op-min} has the form of the equation \eqref{eq:f-min}, then the vector $\bar{c}^* = (c_1^*,\dots,c_n^*)$ satisfies 
	\begin{equation}\label{eq:dif-opt}  
		\nabla_{\bar{c}}\mathcal{M}(\bar{c})= -2n^{-1}\bar{A}K+2\bar{c}K(\lambda I +n^{-1}K)=  0,
	\end{equation}
	for $\bar{c} = \bar{c}^*$. Since the kernel function $\kappa$ is a positive semidefinite function, then the Gram matrix $K$ is a positive semidefinite matrix singular matrix. Thus, we obtain the formula for $\bar{c}^*$:
	\begin{equation}\label{eq:cstar}
		\bar{c}^* = \bar{A}K[K(n\lambda I+ K)]^{\dagger},
	\end{equation}
	where, $[K(n\lambda I+ K)]^{\dagger}$ is the generalized inverse of the matrix $[K(n\lambda I+ K)]$. 
	
	Next, we specify the space $\mathcal{W}$. To ensure the necessary smoothness, we choose the standard Sobolev space $W_2^s(\mathbb{R}^d)$, where $d=4$ and $s=3$. In this case, the Sobolev embedding theorem is satsified for $s > d/2$. This condition implies that $\mathcal{W}$ is a space of continuous functions, and the function values are continuous linear functionals, indicating that it is a reproducing kernel Hilbert space. The reproducing kernel of the Sobolev space $W_2^3(\mathbb{R}^4)$ is given  by \cite{NovUllWozZha18}:
	\begin{equation}\label{eq:kernel-2d}
		\kappa(\mathbf{x}, \mathbf{t}) = \int_{\mathbb{R}^4} \frac{\prod_{j=1}^4 \cos \left(2 \pi\left(x_j - t_j\right) u_j\right) }{(1+\sum_{0<|\alpha|_1 \leq 3} \prod_{j=1}^4(2 \pi u_j)^{2 \mathbf{\alpha}_j})} \mathrm{d} \bf u,
	\end{equation}
	for all $\mathbf{x} = (x_1,\dots,x_4), \mathbf{t} = (t_1,\dots,t_4) \in \mathbb{R}^4$, where  $|\mathbf{\alpha}|_1 = \sum_{j=1}^4\alpha_j$ with non-negative integers $\mathbf{\alpha}_j$.

	
	We have now gathered all the ingredients needed to approximate the scattering amplitude $A$.
	We approximate $A$ by $$f(\mathbf{t}) = \sum^n_{i=1}\bar{c}^*_i \kappa(\mathbf{x}_i,\mathbf{t}),$$ where
	$\mathbf{t} = (t_1,t_2,t_3,t_4)\in\mathbb{R}^4$, $\mathbf{x}_i\in \mathcal{D}$ and $\bar{c}^*$ is given by \eqref{eq:cstar}. 
	The reconstruction of the susceptibility $\eta$ via \eqref{eq:eta_2d} requires differentiation of  $A(\br_1,\br_2)$ with respect to $\br_1$, which in turn requires differentiation of $f (t_1,t_2,t_3,t_4)$ with respect to the first two components  $t_1$ and $t_2$, yielding 
	\begin{equation}\label{Lap2D_f}
		\Laplace_{t_1,t_2} f(\mathbf{t}) =  \sum^n_{i=1}\bar{c}^*_i\kappa^{\prime \prime}(\mathbf{x}_i,\mathbf{t}),
	\end{equation}
	where $\kappa^{\prime \prime}(\mathbf{x},\mathbf{t})$ is given by 
	\begin{equation*}
		\kappa^{\prime \prime}(\mathbf{x}, \mathbf{t}) = -4\pi^2\int_{\mathbb{R}^4} \frac{(u_1^2+u_2^2)\prod_{j=1}^4 \cos \left(2 \pi\left(x_j - t_j\right) u_j\right) }{(1+\sum_{0<|\alpha|_1 \leq 3} \prod_{j=1}^4(2 \pi u_j)^{2 \mathbf{\alpha}_j})} \mathrm{d}\bf u.
	\end{equation*}
	
	Thus, we have identified a functional $f(\mathbf{t})$ within a RKHS $W_2^3(\mathbb{R}^4)$  that approximates the scattering amplitude $A(\br_1,\br_2)$.  Moreover,  the differentiation of $f$ (in turn $A$) is explicitly obtained in terms of \eqref{Lap2D_f}, which is needed for the reconstruction of $\eta$ in \eqref{eq:eta_2d}.
	
	\subsection{Three dimensional problem}\label{subsec:3d} We consider next the case of a three-dimensional medium. The total field $U$ still satisfies the time-harmonic wave equation \eqref{eq:wave} and the integral equation  \eqref{eq:Int_s}, but the three-dimensional Green's function is given by $$\displaystyle G(\br,\br^{\prime}) =\frac{e^{\i k|\br-\br^{\prime}|}}{4\pi|\br-\br^{\prime}|}.$$ 
	In the far field, the Green’s function takes the asymptotic form for $\br \gg \br^{\prime} $ given by,
	\begin{equation}\label{eq:asy-g-3d}
		G\left(\mathbf{r}, \mathbf{r}^{\prime}\right)=\frac{e^{\i k\left|\br\right|}}{4\pi\left|\br\right|}\left( e^{-\i k\hat{\br}\cdot\br^{\prime}} +  \mathcal{O}\left(\frac{1}{\left|\br\right|}\right)\right).
	\end{equation}
	It follow that the asymptotic behavior of the total field is given by,
	\begin{equation}
		\begin{aligned}
			U(\mathbf{r})  &=a\left(\mathbf{r}_1\right) G\left(\mathbf{r}, \mathbf{r}_1\right)+k^2 \int G\left(\mathbf{r}, \mathbf{r}^{\prime}\right) \eta\left(\mathbf{r}^{\prime}\right) U\left(\mathbf{r}^{\prime}\right) d \mathbf{r}^{\prime} \\
			& \sim \frac{e^{\i k\left|\br\right|}}{4\pi\left|\br\right|} e^{-\mathrm{i} k \hat{\mathbf{r}} \cdot \mathbf{r}_1} a\left(\mathbf{r}_1\right)+\frac{e^{\i k\left|\br\right|}}{4\pi\left|\br\right|} k^{2} \int e^{-\mathrm{i} k \hat{\mathbf{r}} \cdot \mathbf{r}^{\prime}} \eta\left(\mathbf{r}^{\prime}\right) U\left(\mathbf{r}^{\prime}\right) d \mathbf{r}^{\prime}
		\end{aligned}
	\end{equation}
	In the far-zone of the scatterer, the field behaves as an outgoing spherical wave:
	\begin{equation}\label{eq:asy-u-3d}
		U(\mathbf{r}) \sim \frac{e^{\i k\left|\br\right|}}{4\pi\left|\br\right|} \left(a(\br_1) e^{-\mathrm{i} k \hat{\mathbf{r} \cdot \mathbf{r}_1}}+A(\mathbf{r})\right),
	\end{equation}
	where the scattering amplitude $A$ is defined by
	\begin{equation}\label{eq:def-ScatAmp-3d}
		A(\mathbf{r})=k^{2} \int_{\Omega} e^{-\i k \hat{\mathbf{r}} \cdot \mathbf{r}^{\prime}} \eta\left(\mathbf{r}^{\prime}\right) U\left(\mathbf{r}^{\prime}\right) d \mathbf{r}^{\prime}.
	\end{equation}
	
	Following the same procedure as before, the scattering amplitude in three dimensions is rewritten as,
	\begin{equation}\label{eq:ScatAmp-r1r2-3d}
		\begin{aligned}
			A(\br_1,\br_2) &= a(\br_1) k^{2}  \int_{\Omega} e^{-\i k \widehat{\br}_2 \cdot \br} \eta(\br) G(\br , \br_1 ) d\br,
		\end{aligned}
	\end{equation}
	where we assume the sources have unit amplitude. The inversion formula for $\eta$ is thus given by,
	\begin{align}\label{eq:eta_3d}
		\eta(\br_1) = -\frac{e^{\i k \widehat{\br}_2 \cdot \br_1}}{k^{2   }} \left(   \Laplace_{\br_1} A(\br_1,\br_2) + k^2 A(\br_1,\br_2) \right).
	\end{align}
	
	The three-dimensional inverse problems is also considered in \cite{GilLevScho18}, where \eqref{eq:ScatAmp-r1r2-3d} is converted into a linear system and then a regularized pseudoinverse solution is used to recover the susceptibility. In addition, the number of sources, detectors, and volume elements in $\Omega$ was chosen so that the resulting linear system is overdetermined. Moreover, the two-dimensional case was not considered in \cite{GilLevScho18}.    
	
	The forward problem of solving \eqref{eq:def-ScatAmp-3d} was carried out by using the coupled-dipole method
	\cite{LevMar16, PurPen73}.
	We discretized the  three-dimensional domain $\Omega$ into voxels $C_k\left(k=1,2, \cdots N_\nu\right) $ of volume  $h^3$, where $N_\nu $ is the total number of voxels. We split the integral of  \eqref{eq:ScatAmp-r1r2-3d} into singular- and non-singular parts according to
	\begin{align}\label{A_int3D}
		k^{2} \int_{\Omega} e^{-\i k \widehat{\br}_2 \cdot \br} \eta(\br) G(\br , \br_1 ) d\br=k^{2} \sum_{k=1}^{N_\nu} \int_{\Omega_k \backslash \Omega_s} e^{-\i k \widehat{\br}_2 \cdot \br} \eta(\br) G(\br , \br_1 ) d\br+k^{2} \int_{\Omega_s} e^{-\i k \widehat{\br}_2 \cdot \br} \eta(\br) G(\br , \br_1 ) d\br ,
	\end{align}
	where we assume that $k h \ll 1$. Combining the fact that $e^{\i k\left| \br-\br_1\right|} \approx 1 + \i k \left| \br-\br_1\right|$, the singular part of \eqref{A_int3D} is computed as
	\begin{align*}
		\begin{aligned}
			\int_{\Omega_s} e^{-\mathrm{i} k \widehat{\mathbf{r}}_2 \cdot \mathbf{r}} \eta(\mathbf{r})   G\left(\mathbf{r}, \mathbf{r}_1\right) d \mathbf{r} &\approx  \frac{1}{4\pi}\eta(\widetilde{\mathbf{r}}) e^{-\mathrm{i} k \widehat{\mathbf{r}}_2 \cdot \widetilde{\mathbf{r}}} \int_{-\frac{h}{2}}^{\frac{h}{2}}\int_{-\frac{h}{2}}^{\frac{h}{2}} \int_{-\frac{h}{2}}^{\frac{h}{2}}\left( 
			\frac{1}{\left| \br-\br_1\right|} +\i k       \right) d \mathbf{r} 
			=\frac{1}{4\pi}\eta(\widetilde{\mathbf{r}}) e^{-\mathrm{i} k \widehat{\mathbf{r}}_2 \cdot \widetilde{\mathbf{r}}}\zeta,  
		\end{aligned}
	\end{align*}
	where 
	$ \zeta= h^2(\xi +\i k h)$ with $ \xi = \text{ln}(26+15\sqrt{3})-\frac{\pi}{2} \approx 2.38,$ and $\widetilde{\mathbf{r}}$ is the central point in $\Omega_s$.
	
	Following along similar lines as in the previous subsection, 
	we reconstruct the three-dimensional susceptibility $\eta$ based on \eqref{eq:eta_3d}. We consider the optimization problem \eqref{eq:op-min} in the Sobolev space 
	$W_2^s(\mathbb{R}^d)$, where $d=6$ and $s=4$, 
	which is a RKHS. 
	The space $W_2^4(\mathbb{R}^6)$ has the 
	reproducing kernel given  by \cite{NovUllWozZha18}:
	\begin{equation}\label{eq:kernel-3d}
		\kappa(\mathbf{x}, \mathbf{t}) = \int_{\mathbb{R}^6} \frac{\prod_{j=1}^6 \cos \left(2 \pi\left(x_j - t_j\right) u_j\right) }{(1+\sum_{0<|\alpha|_1 \leq 4} \prod_{j=1}^6(2 \pi u_j)^{2 \mathbf{\alpha}_j})} \mathrm{d}\bf u,
	\end{equation}
	where $\mathbf{x} = (\mathbf{r}_1,\mathbf{r}_2) = \mathcal{D} \in \mathbb{R}^6$ with $\mathbf{r}_1$ and $\mathbf{r}_2$ are the source and detector positions, respectively. 
	Similarly, we identify a functional $f$ within the RKHS $W_2^4(\mathbb{R}^6)$  that approximates the scattering amplitude  $A$.
	The functional  $$f(\mathbf{t}) = \sum^n_{i=1}\bar{c}^*_i\kappa(\mathbf{x}_i,\mathbf{t}),$$ where $\mathbf{t} = (t_1,t_2,t_3,t_4,t_5,t_6)\in\mathbb{R}^6$, $\mathbf{x}_i\in \mathcal{D}$ and $\bar{c}^*$ is as in  \eqref{eq:cstar}. 
	The reconstruction of $\eta$ from \eqref{eq:eta_3d} requires differentiation of  $A(\br_1,\br_2)$ with respect to $\br_1$, which in turn requires differentiation of $f (t_1,t_2,t_3,t_4,t_5,t_6)$ with respect to $t_1$, $t_2$ and $t_3$, yielding 
	\begin{equation}\label{Lap3D_f}
		\Laplace_{t_1,t_2,t_3} f(\mathbf{t}) =  \sum^n_{i=1}\bar{c}^*_i \kappa^{\prime \prime}(\mathbf{x}_i,\mathbf{t}),
	\end{equation}
	where $\kappa^{\prime \prime}(\mathbf{x},\mathbf{t})$ is given by 
	\begin{align*}
		\displaystyle 
		\kappa^{\prime \prime}(\mathbf{x}, \mathbf{t}) = -4\pi^2\int_{\mathbb{R}^6} \frac{(u_1^2+u_2^2+u_3^2)\prod_{j=1}^6 \cos \left(2 \pi\left(x_j - t_j\right) u_j\right) }{(1+\sum_{0<|\alpha|_1 \leq 4} \prod_{j=1}^6(2 \pi u_j)^{2 \mathbf{\alpha}_j})} \mathrm{d}\bf u.
	\end{align*}
	Finally, following the same procedure as in the previous subsection, we reconstruct $\eta$ using \eqref{eq:eta_3d}.
	
	%
	
	For three-dimensional reconstructions, we present two methods. First we perform two-dimensional reconstructions using \eqref{eq:eta_2d} for each slice, and then make a three-dimensional reconstruction by assembling these cross-sectional reconstructions together. The second, 
	is to  directly use \eqref{eq:eta_3d} for three-dimensional reconstruction.

	\section{numerical experiments}\label{sec:numerics}
	
	To illustrate the reconstruction method, we present three numerical experiments.
	
	\subsection{Reconstructions of a three-ball model}  We consider reconstructing a three-ball model \cite{GilLevScho18}. The model system consists of a volume of dimensions
	70 nm $\times$ 70 nm $\times$ 40 nm. Three spherical scatterers of radius 12 nm are placed so their lowest point is 3 nm from the bottom of the sample. The minimum distance between any two spheres is 5 nm. The susceptibility of two of the spheres is set to 1.275, while the remaining sphere has a susceptibility of 1.885.
	The wave number $k = 2\pi/500 = 0.0126$.
	Gaussian noise at the 1$\%$ level is added to the scattering amplitude and
	the positions of the sources. We divide the system into 17 layers, thus each layer is a two-dimensional medium, we solve the forward problem on each layer with a discretization of  35$\times$35, while the susceptibility for each layer is reconstructed on a grid of size 31$\times$31 by \eqref{eq:eta_2d}. This model system is similar to what is considered in \cite{GilLevScho18}. We use 150 randomly placed sources and 7 detectors of $L$ shape in this experiment, for each layer, to show the feasibility of the  two-dimensional reconstruction method.
	
	\begin{figure}[!htb]
		\centering
		\ \
		\subfigure[ \tiny{ Positions of 150 sources}]{\includegraphics[width=0.215\textwidth]{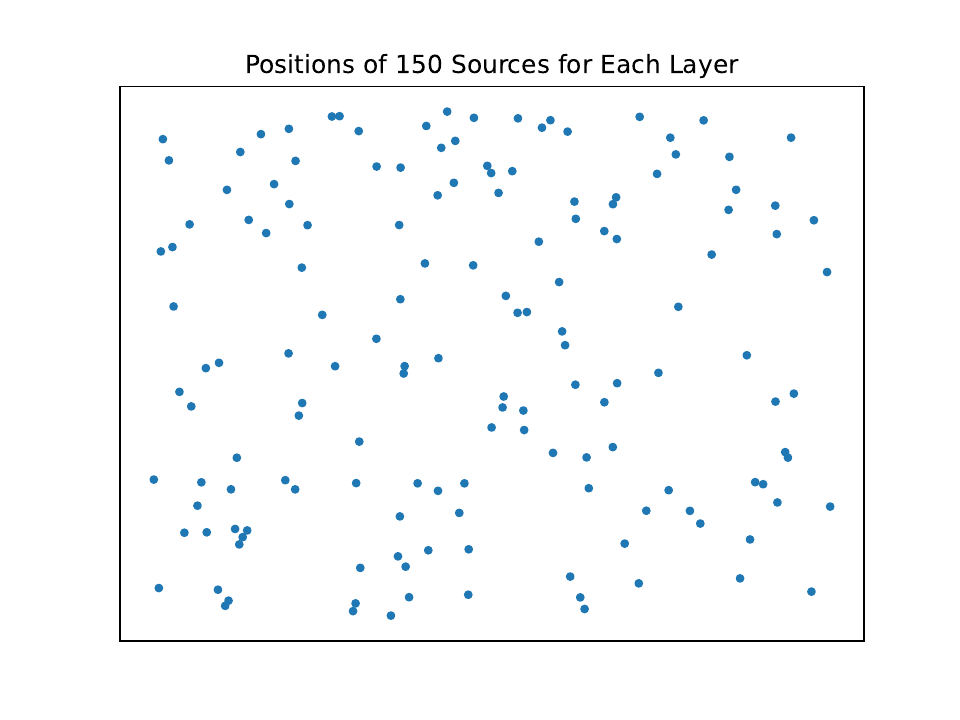}}  \ \  \ \ \ \ \
		\subfigure[ \tiny{ layer 1-5 and 13-17 }]{\includegraphics[width=0.302\textwidth]{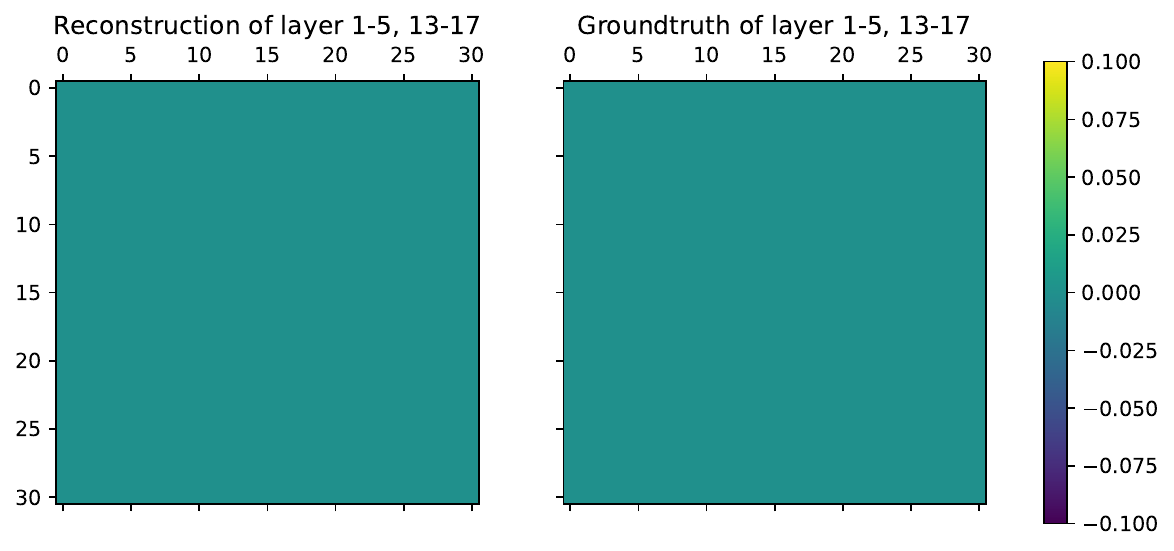}}
		\subfigure[ \tiny{ layer 6 and 12 }]{\includegraphics[width=0.302\textwidth]{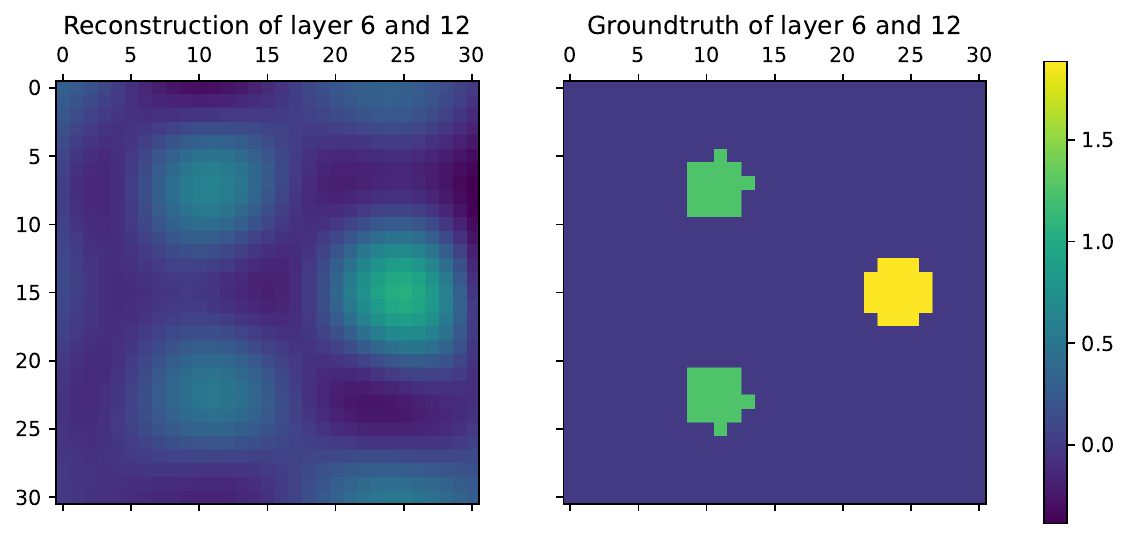}}\\

		\subfigure[ \tiny{ layer 7 and 11 }]{\includegraphics[width=0.302\textwidth]{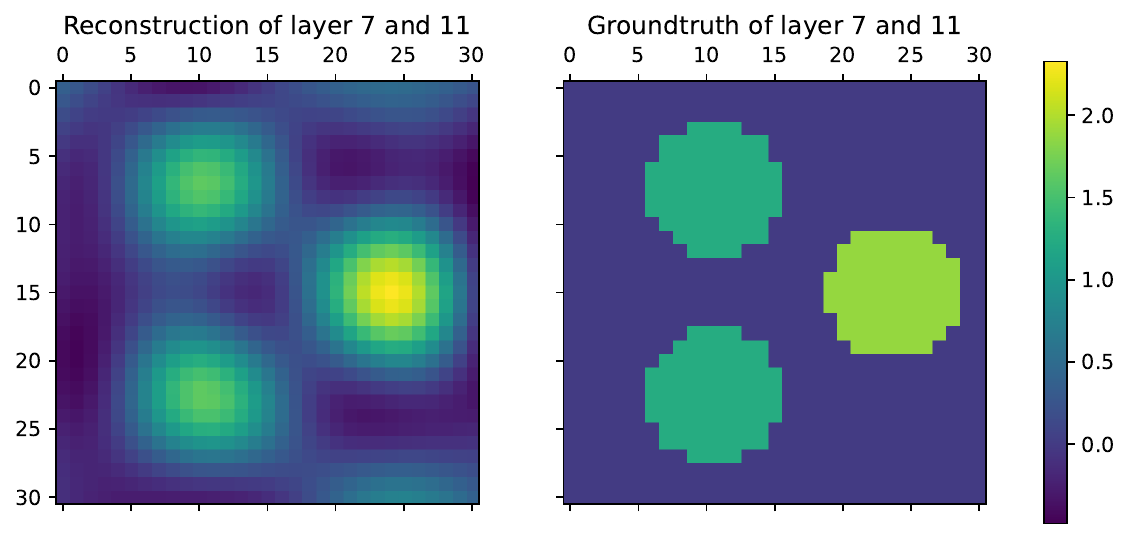}}
		\subfigure[ \tiny{ layer 8 and 10 }]{\includegraphics[width=0.302\textwidth]{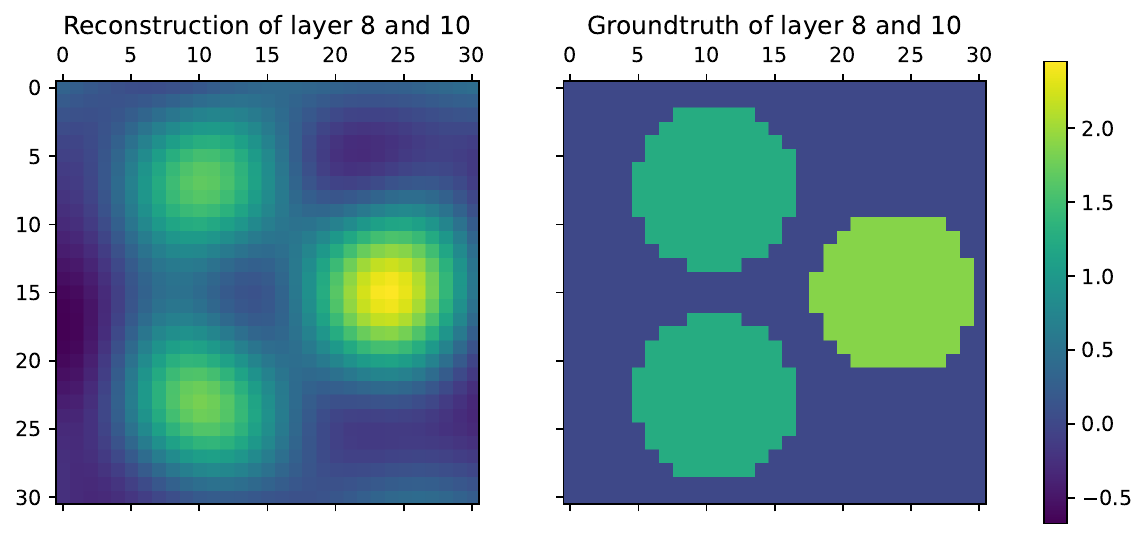}} 
		\subfigure[ \tiny{ layer 9 }]{\includegraphics[width=0.302\textwidth]{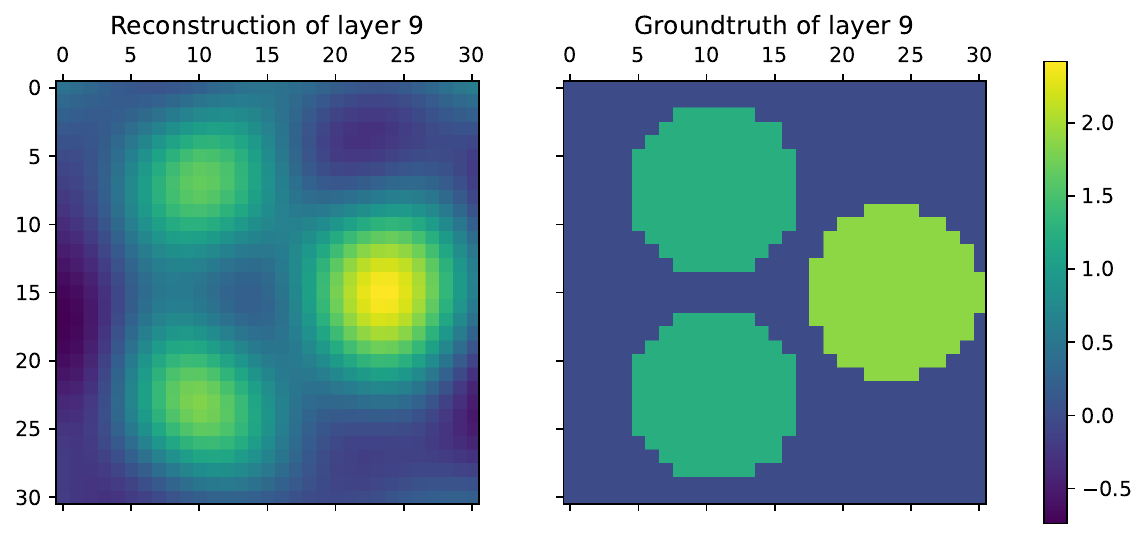}}

		\caption{\small \it Susceptibility reconstruction of three spherical scatterers on each layer by the two-dimensional method.
		} \label{fig:150-7-2d-slice}
	\end{figure}
	
	Figure \ref{fig:150-7-2d-slice} shows the position of the sources and the reconstructions for each layer. 
	Due to the symmetry of the medium and the source positions for each layer being the same, some layer reconstructions are also identical. In our three-ball model,  layer 9 is the central layer, so that layer 8 $\&$ 10, and layer 7 $\&$ 11 (see Figure \ref{fig:150-7-2d-slice} (d) and (e) ) yield identical layer reconstructions. Thus, one can reconstruct on either of the two layers, hence reducing the reconstruction time. It is also worth noting that the Gram matrix $K$ for calculating the vector $\bar{c}^*$ in \eqref{eq:cstar} for each layer is the same, since the position of the sources and the reconstruction grid are the same for each layer. This also applies to subsection \ref{subsec:2d-neuron}.
	
	\subsection{Reconstructions of a neuron} We consider a neuron from a rat brain, of size 3.928$\mu m$$\times$5.44$\mu m$ $\times$2.248$\mu m$. The image was obtained from experiment and manually segmented \cite{MarDeeYamBusEll00,MarGupWonQiaSos02}. 
	The main trunk of the neuron has a susceptibility of $\eta=0.25$ and the branches spines have a susceptibility of $\eta=0.38$, while the background has susceptibility  $\eta=0$. Once again, 1$\%$ Gaussian noise is added to the scattering amplitude and the positions of the sources. Both two-dimensional method and three-dimensional methods were tested. 
	
	\subsubsection{Two-dimensional reconstruction method}\label{subsec:2d-neuron}
	We manually segmented the neuron into a grid of size 33$\times$46$\times$18, then we calculated the forward problem on 18 layers with each layer of size 33$\times$46. Subsequently, we reconstructed each of these 18 layers onto a 31$\times$44 grid, denoted as  $\mathcal{G}:=\{(g_i,g_j)\}$, with $i=1,\dots, 31$, and $j =1,\dots,44$. For each layer, 500 randomly spaced sources and 7 detectors were employed. 
	
\begin{figure*}
		\centering
		\subfigure[ \tiny{ Positions of 500 sources}]{\includegraphics[width=0.15\textwidth]{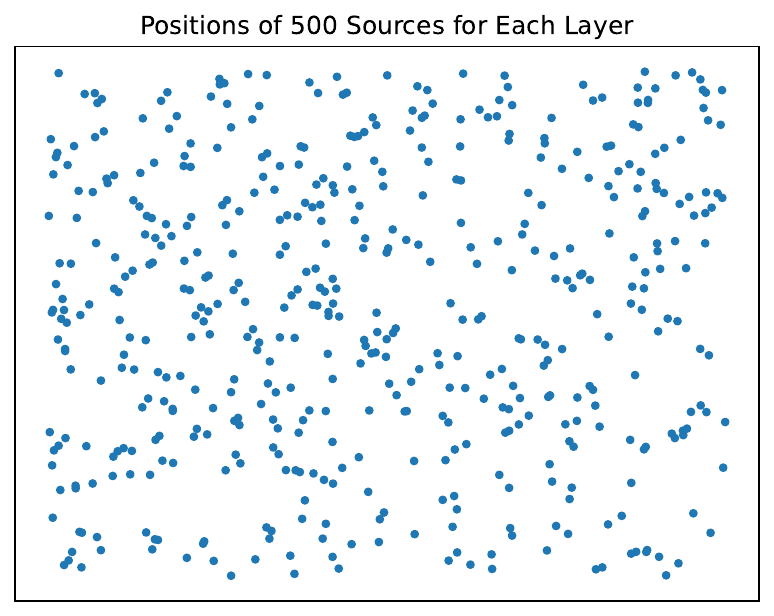}}  
		\subfigure[ \tiny{ layer 1 and 18}]{\includegraphics[width=0.15\textwidth]{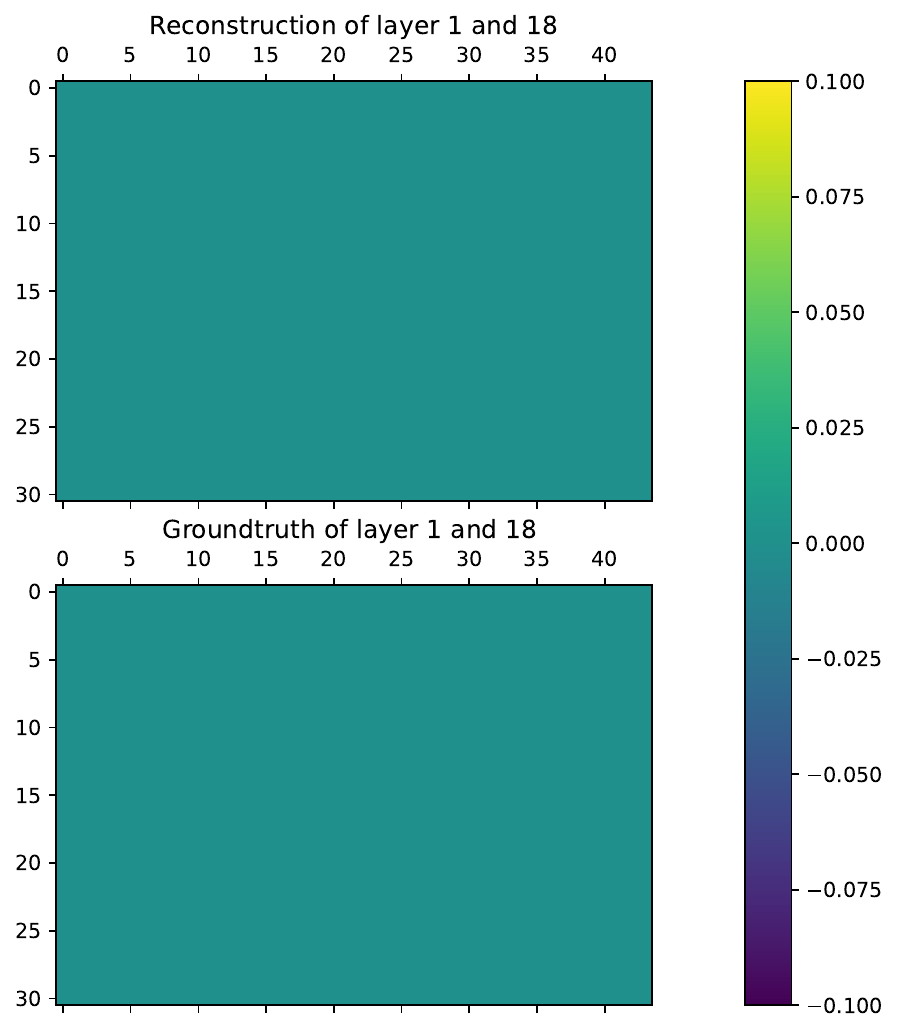}}
		\subfigure[ \tiny{ layer 2  }]{\includegraphics[width=0.15\textwidth]{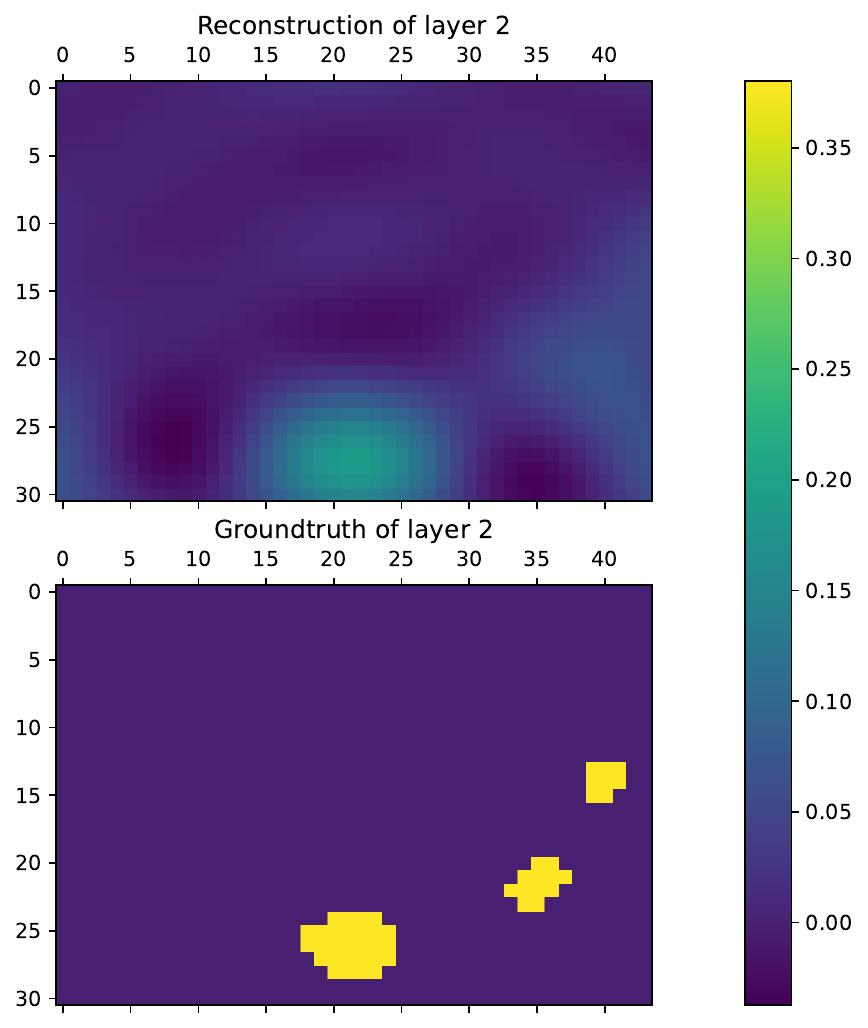}}
		\subfigure[ \tiny{ layer 3  }]{\includegraphics[width=0.15\textwidth]{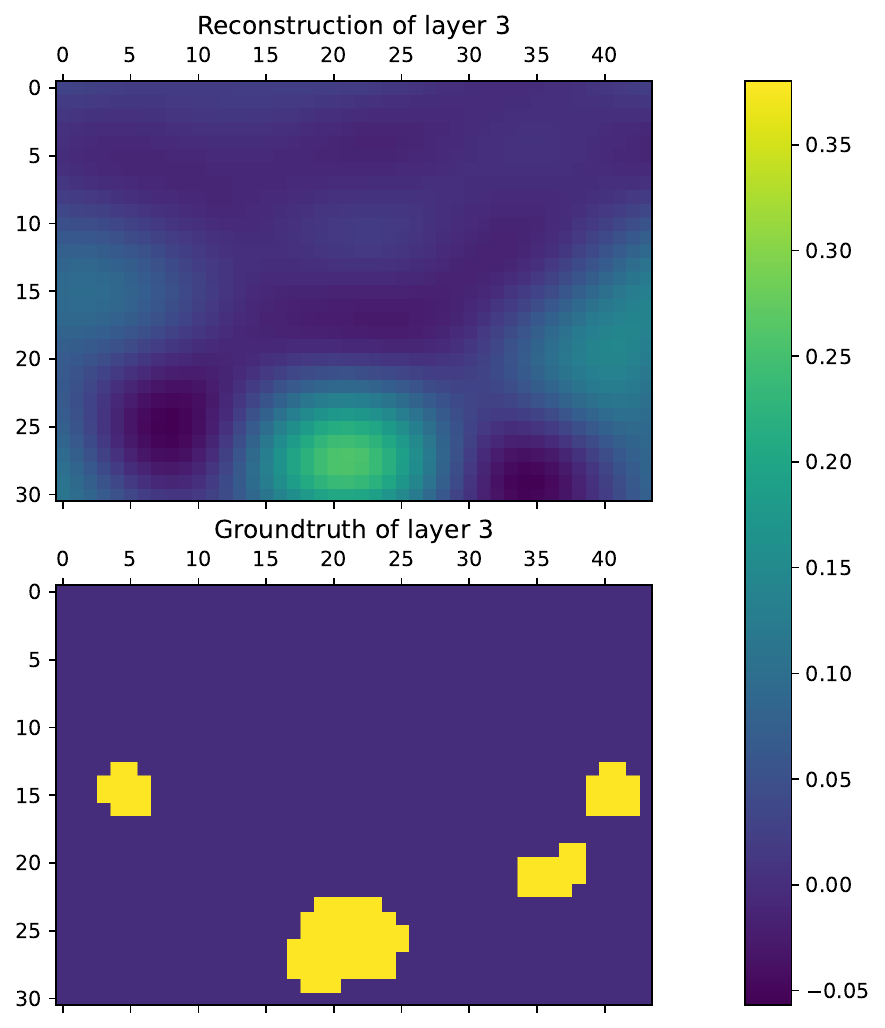}}
		\subfigure[ \tiny{ layer 4  }]{\includegraphics[width=0.15\textwidth]{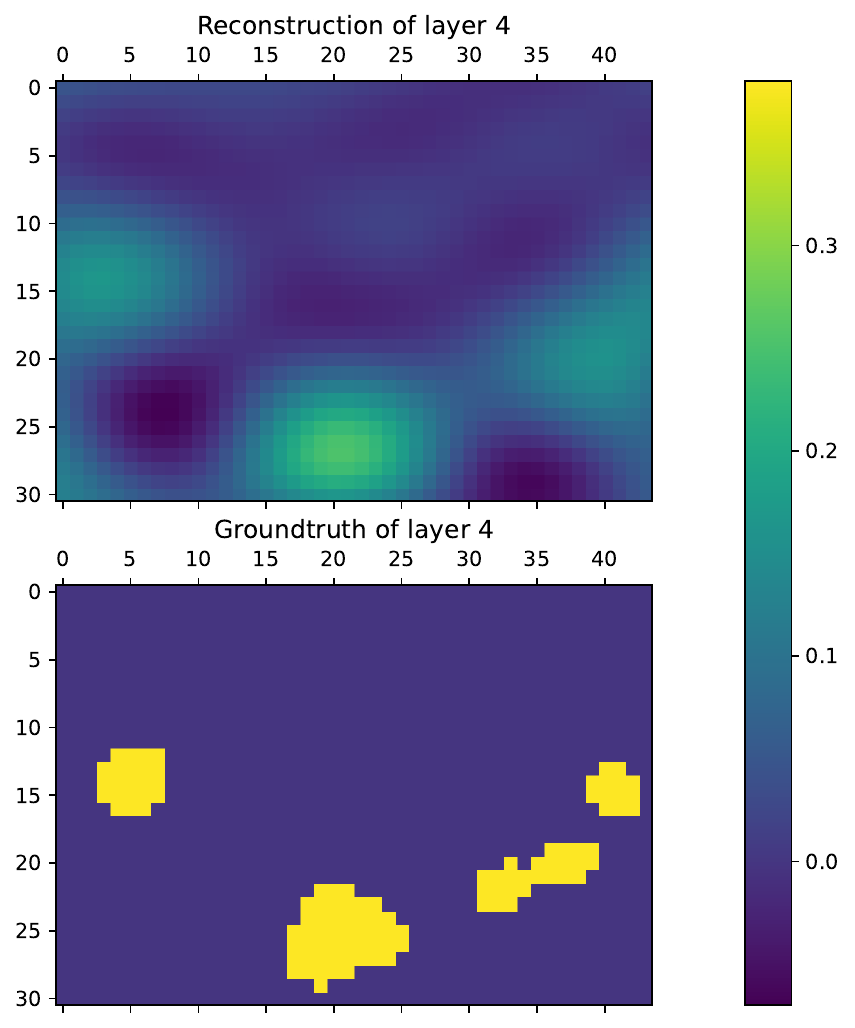}}
		\subfigure[ \tiny{ layer 5  }]{\includegraphics[width=0.15\textwidth]{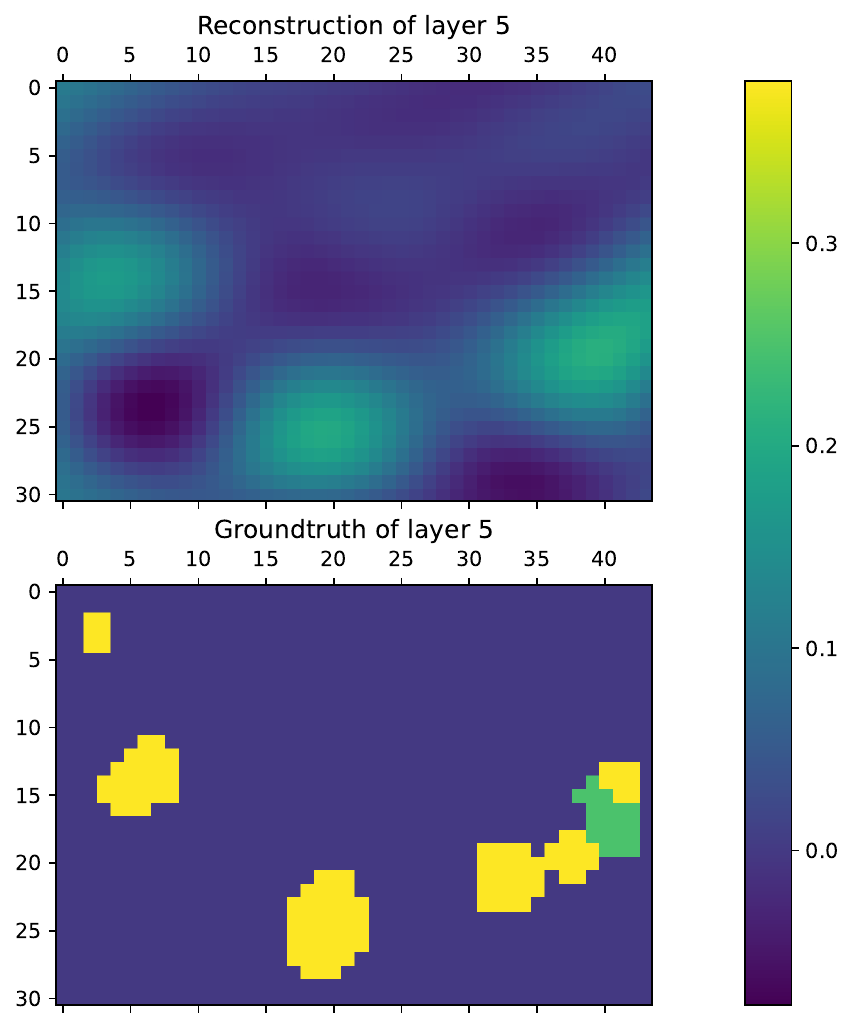}}\\
		
		\subfigure[ \tiny{ layer 6}]{\includegraphics[width=0.15\textwidth]{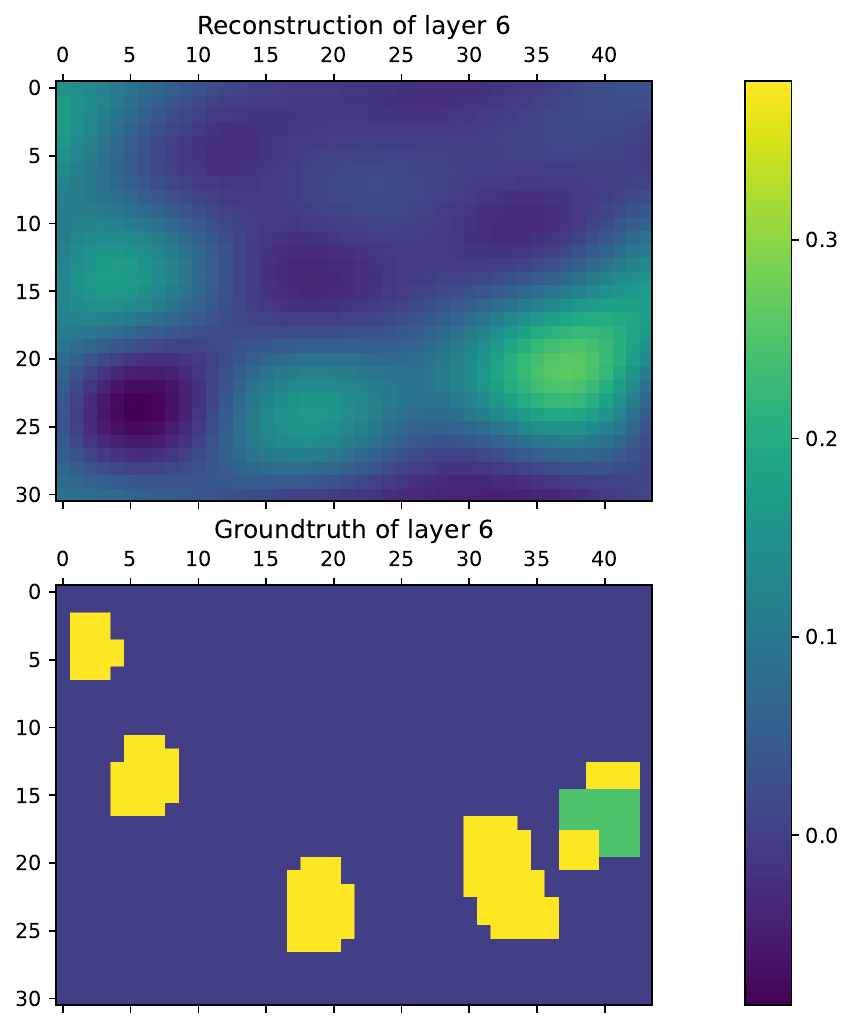}}  
		\subfigure[ \tiny{ layer 7}]{\includegraphics[width=0.15\textwidth]{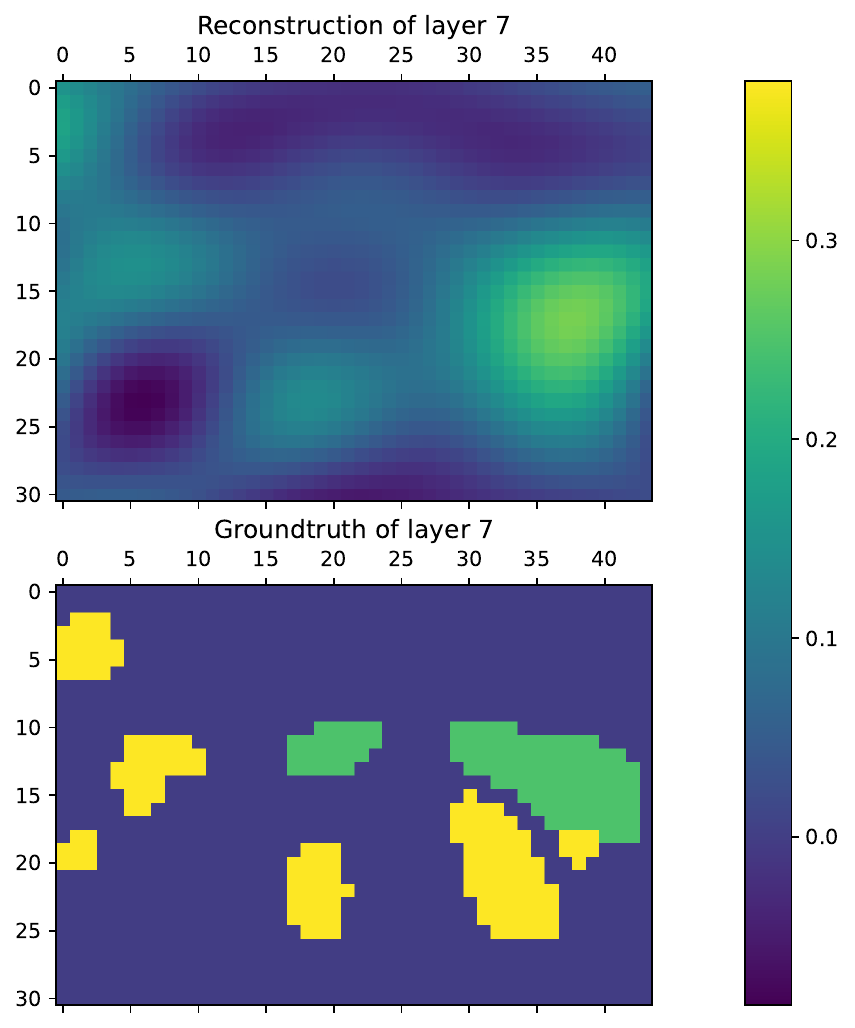}}
		\subfigure[ \tiny{ layer 8  }]{\includegraphics[width=0.15\textwidth]{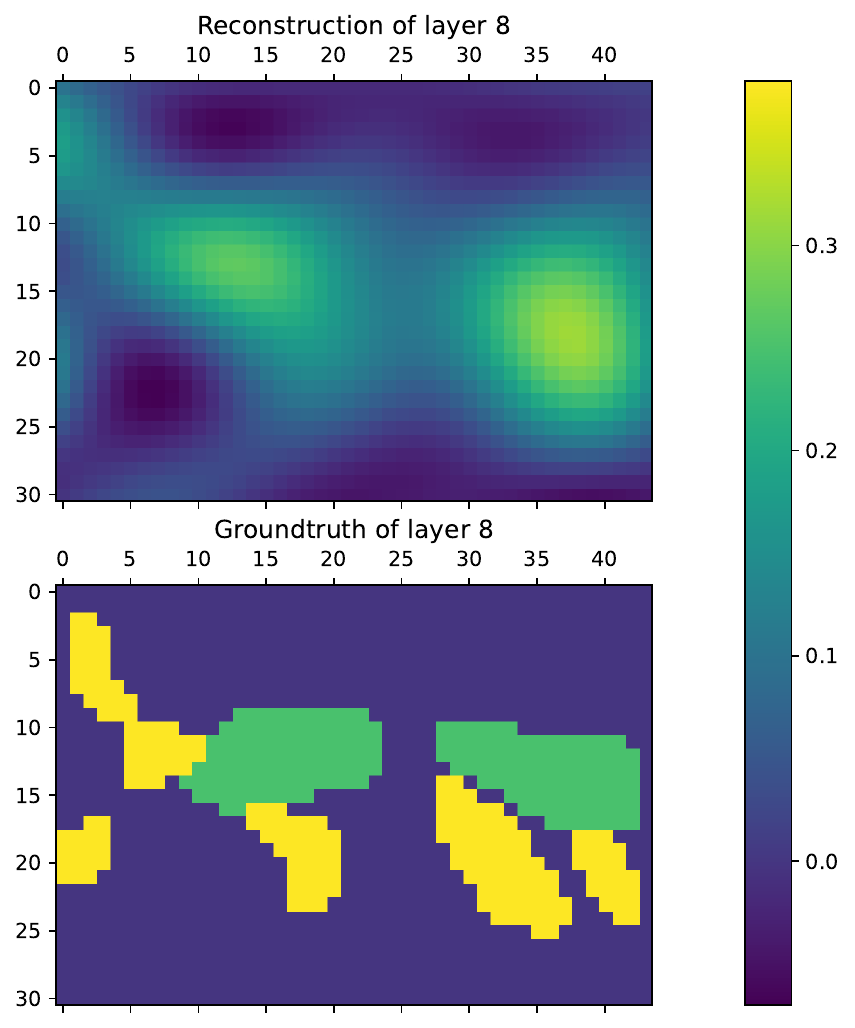}}
		\subfigure[ \tiny{ layer 9  }]{\includegraphics[width=0.15\textwidth]{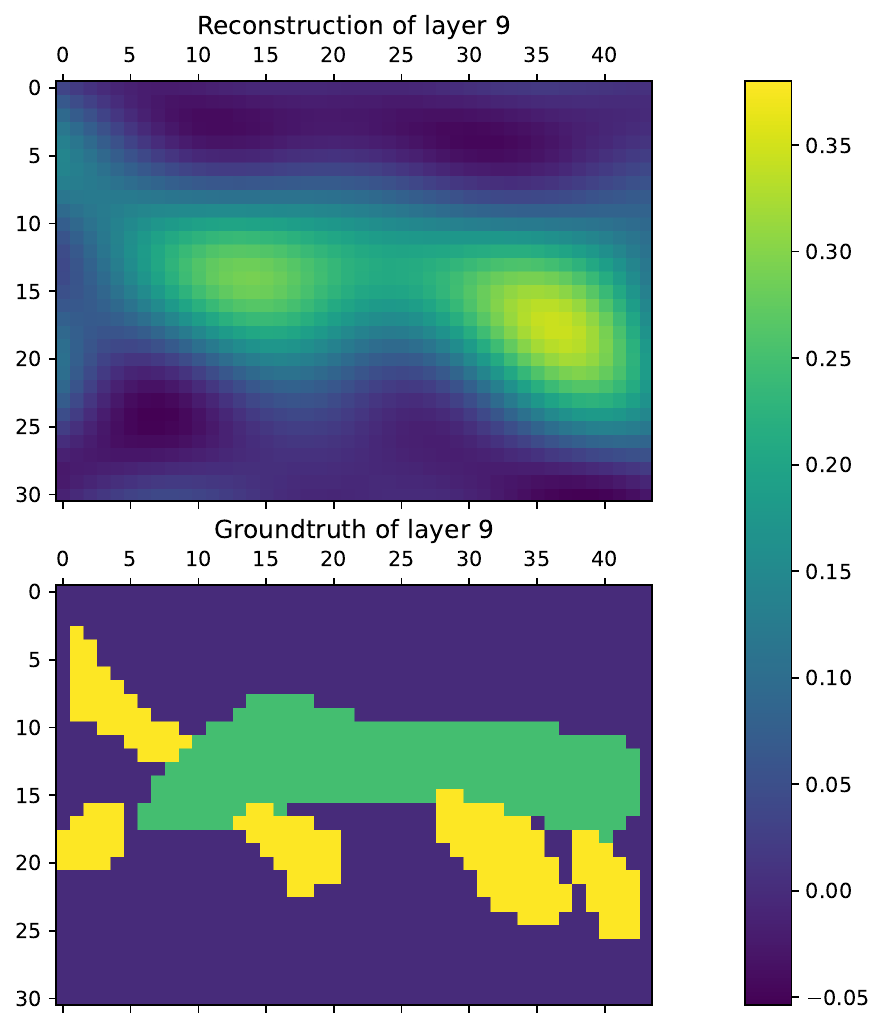}}
		\subfigure[ \tiny{ layer 10  }]{\includegraphics[width=0.15\textwidth]{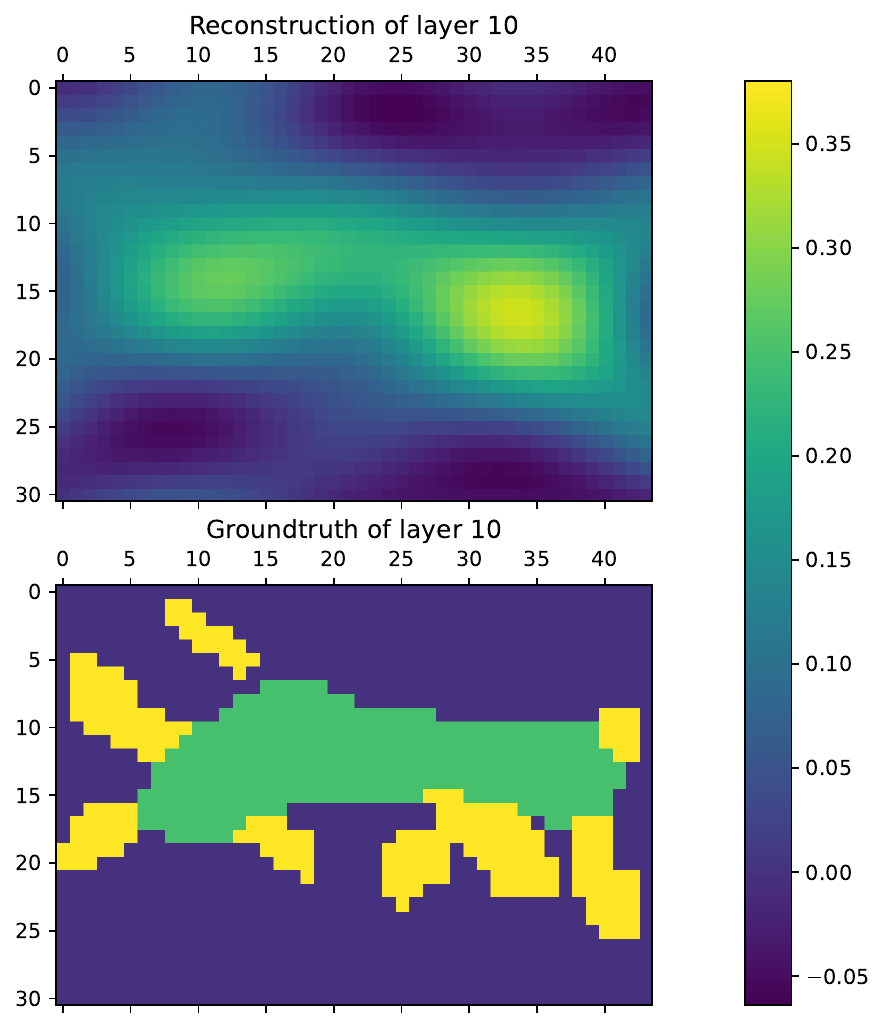}}
		\subfigure[ \tiny{ layer 11  }]{\includegraphics[width=0.15\textwidth]{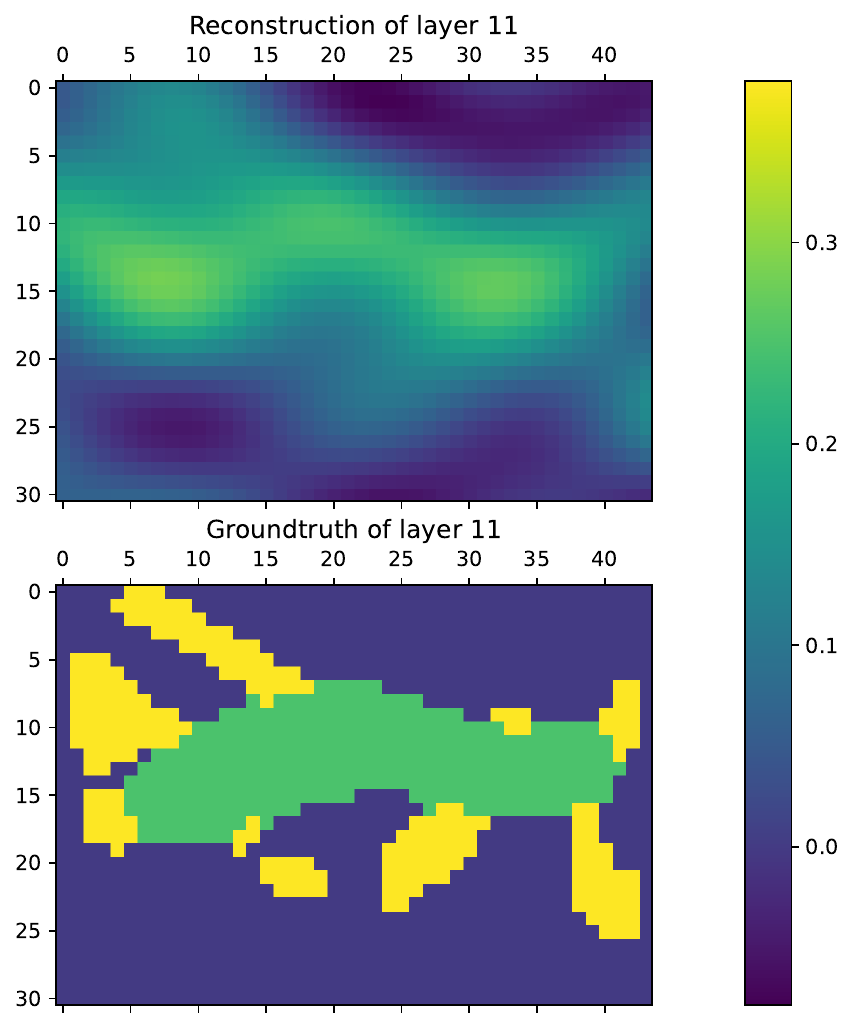}}\\
		\subfigure[ \tiny{ layer 12}]{\includegraphics[width=0.15\textwidth]{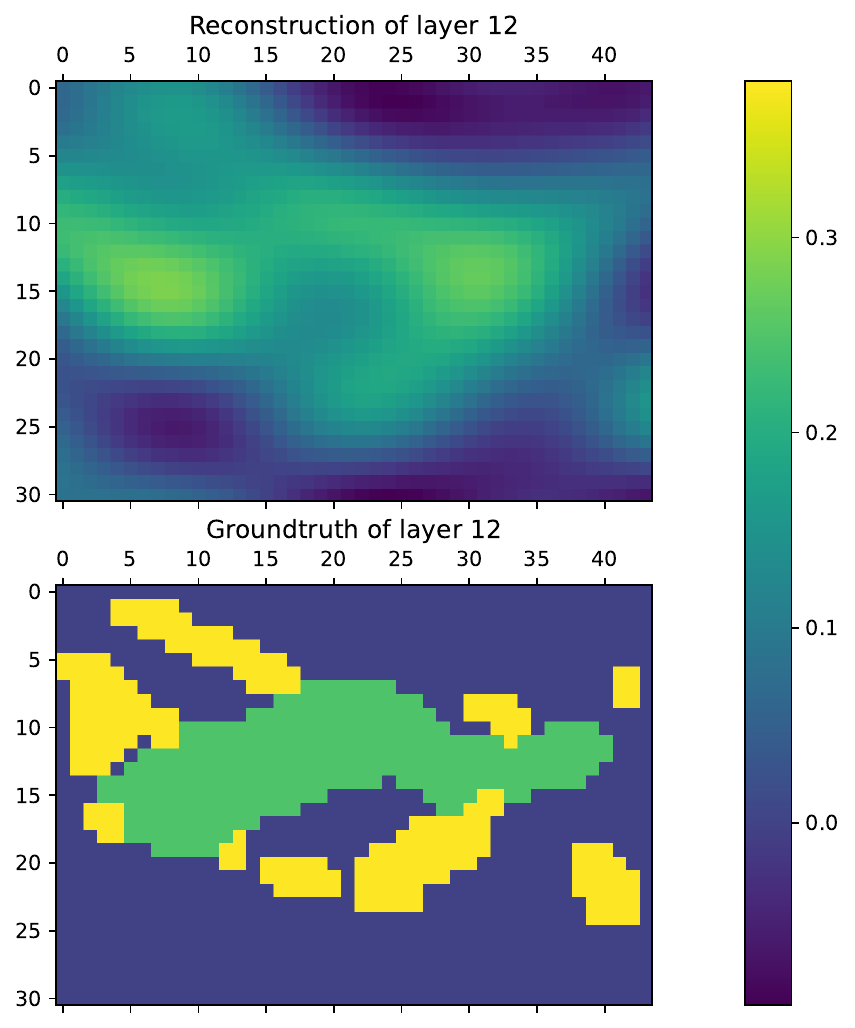}}  
		\subfigure[ \tiny{ layer 13}]{\includegraphics[width=0.15\textwidth]{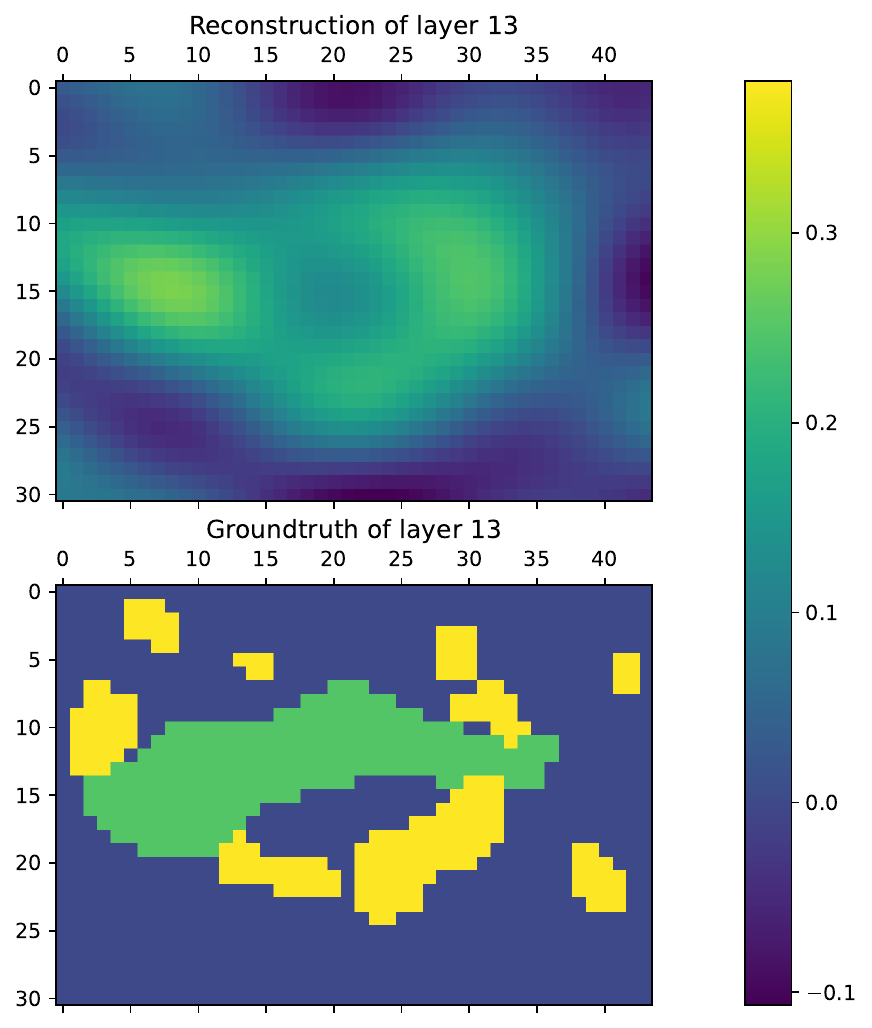}}
		\subfigure[ \tiny{ layer 14  }]{\includegraphics[width=0.15\textwidth]{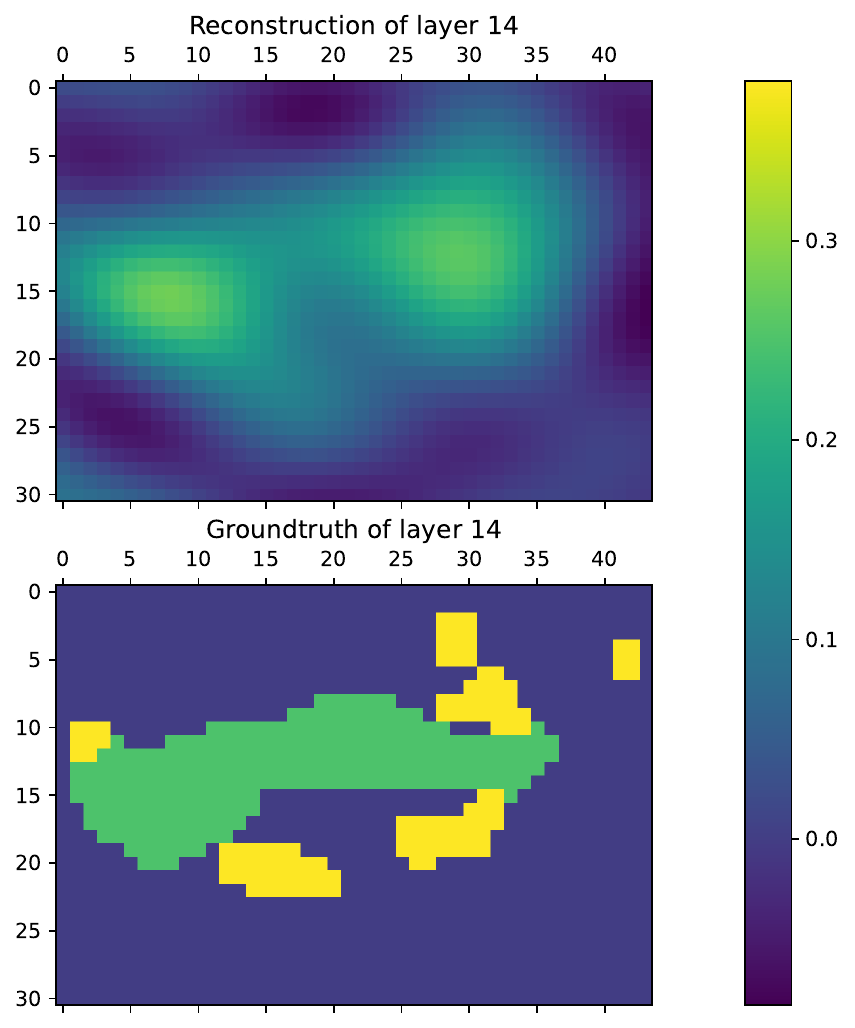}}
		\subfigure[ \tiny{ layer 15  }]{\includegraphics[width=0.15\textwidth]{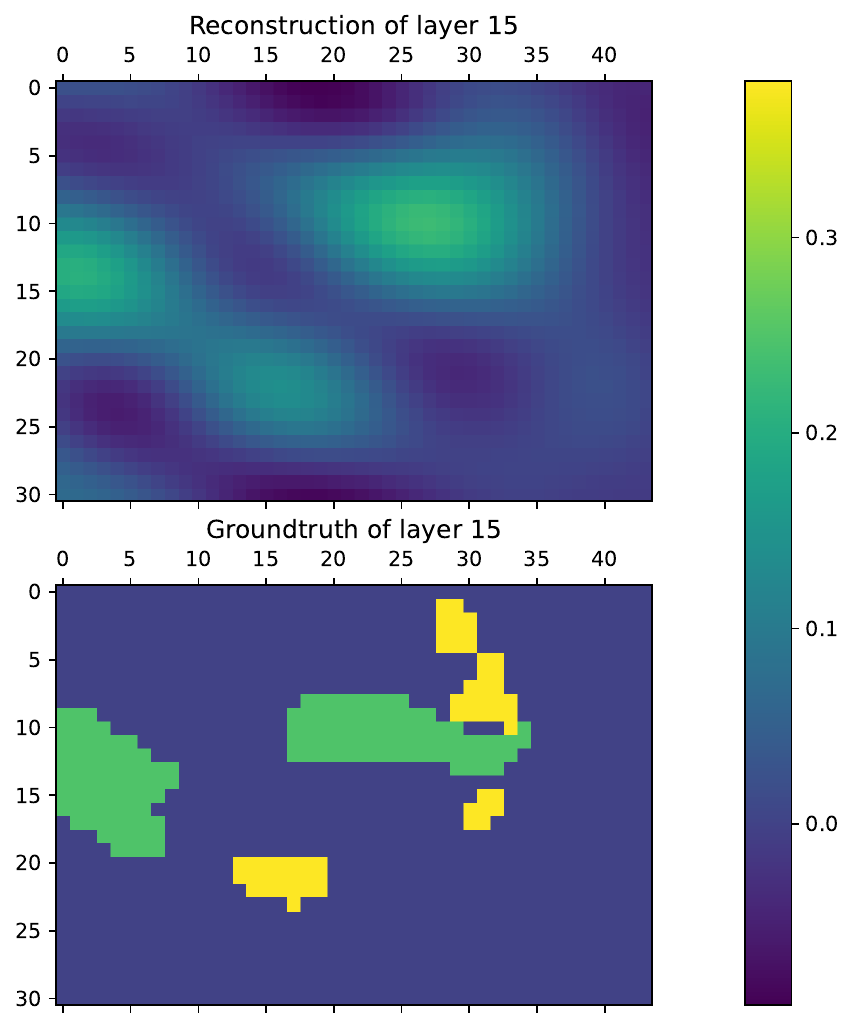}}
		\subfigure[ \tiny{ layer 16  }]{\includegraphics[width=0.155\textwidth]{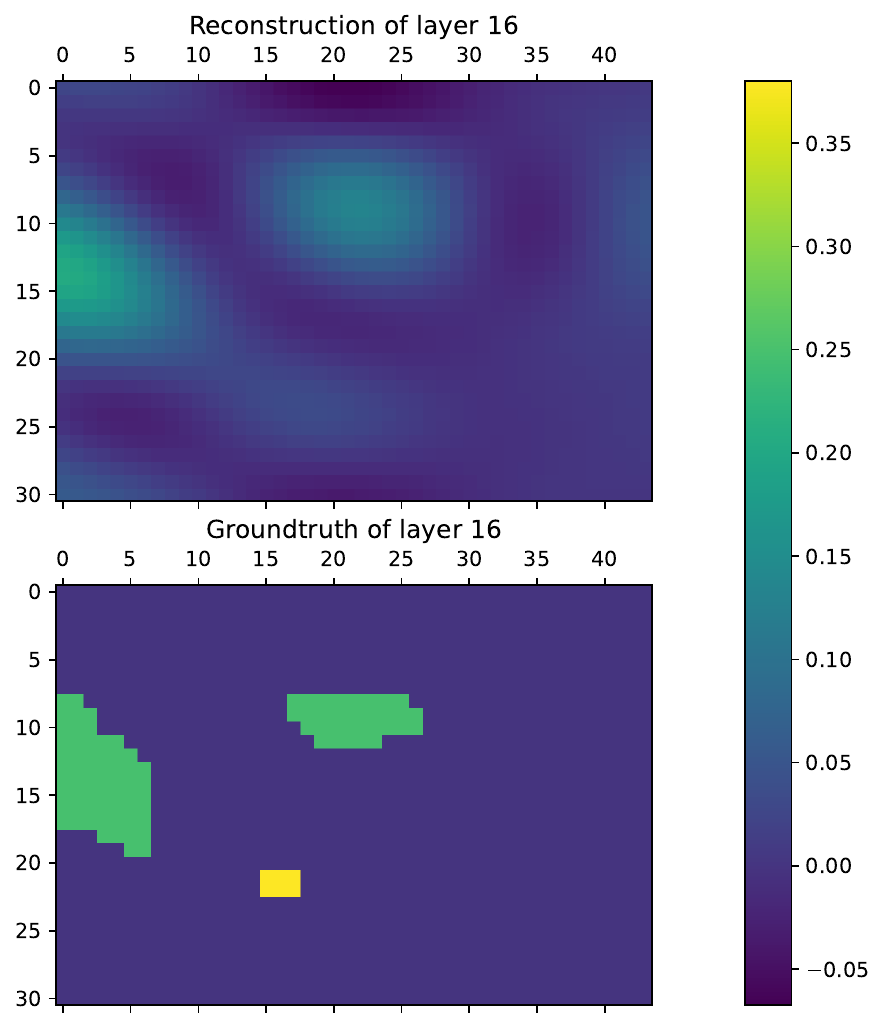}}
		\subfigure[ \tiny{ layer 17  }]{\includegraphics[width=0.15\textwidth]{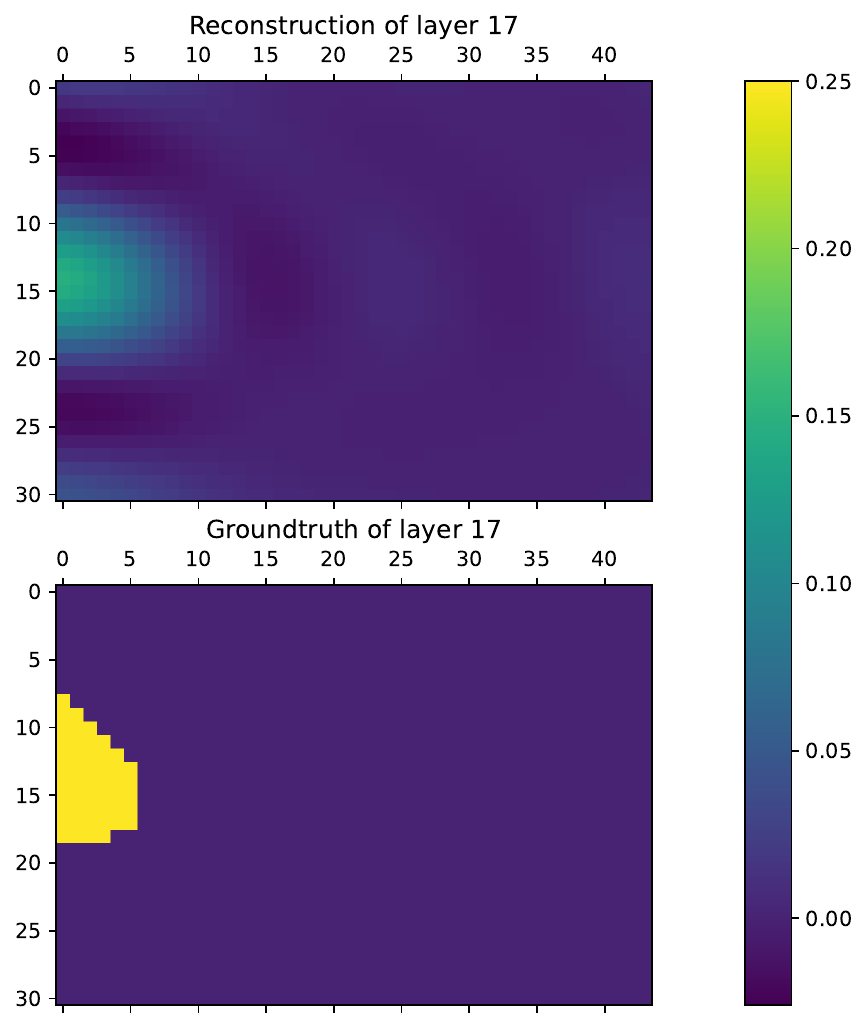}}\\

		\caption{\small \it Neuron model and susceptibility reconstructions on each layer by the two-dimensional method.} \label{fig:500-7-2d-slice}
	\end{figure*}

	Since we approximate the scattering amplitude by $f(\mathbf{t}) = \sum^n_{i=1}\bar{c}^*_i\kappa(\mathbf{x}_i,\mathbf{t})$, $\mathbf{t} = (t_1,t_2,t_3,t_4)\in\mathbb{R}^4$, $\mathbf{x}_i\in \mathcal{D}$, we can approximate the value at grid points by letting $(t_1,t_2) = (g_i,g_j)$ and randomly set $(t_3,t_4)$. In Figure \ref{fig:500-7-2d-slice}, we only set one $(t_3,t_4)$, which is the coordinate of one detector, yielding one susceptibility reconstruction for each grid point.
	
\begin{figure}[!htb]
		\centering
		{\includegraphics[width=0.415\textwidth]{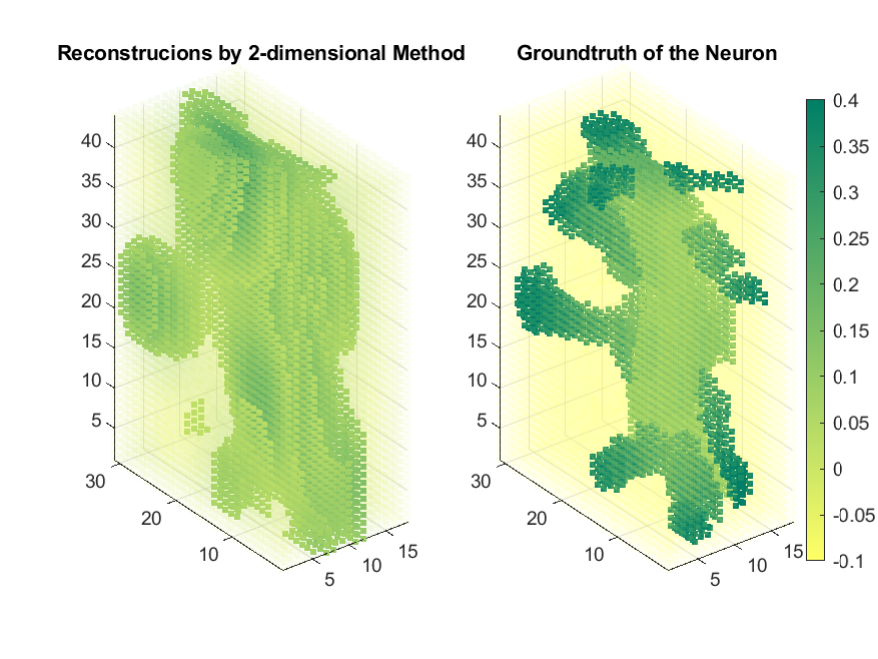}}  
		\caption{\scriptsize \it Neuron reconstructions by two-dimensional method compared to the groundtruth.} \label{fig:500-7-2d-whole}
	\end{figure}

	Then putting all of the layers together, we obtain the reconstruction shown in Figure \ref{fig:500-7-2d-whole}. To better display the three dimensional colour map, we add transparency to the points of value less than 0.085.

	Let us define the relative error for the scattering amplitude:
	$$
	\chi^2=\frac{\sum_{{\br}_1, {\br}_2}\left|A\left({\br}_1, {\br}_2\right)-A_{\mathrm{rec}}\left({\br}_1, {\br}_2\right)\right|^2}{\sum_{{\br}_1, {\br}_2}\left|A\left({\br}_1, {\br}_2\right)\right|^2},
	$$ 
	where $A_{\mathrm{rec}} $ is the reconstructed scattering amplitude. We define the mean relative error of the susceptibility at grid $\mathcal{G}$ as 
	$$
	\delta = \frac{\sum^{31}_{i=1}\sum^{44}_{j=1}|\eta\{g_i,g_j\}-\eta_{rec} \{g_i,g_j\} |}{ 31\times44   },
	$$
	where, $\eta\{g_i,g_j\}$ is the groundtruth susceptibility at grid point $(g_i,g_j)$, and $\eta_{rec}$ is the reconstructed susceptibility. 
	
\begin{figure}[!t]
		\centering
		\subfigure{\includegraphics[width=0.3515\textwidth]{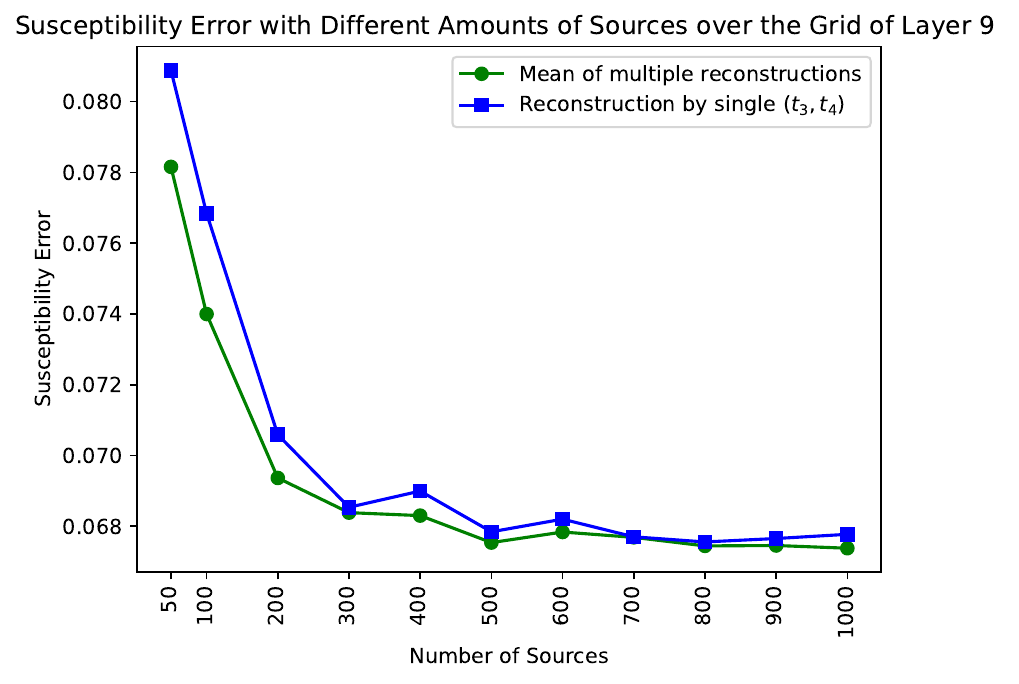}}
		\subfigure{\includegraphics[width=0.3515\textwidth]{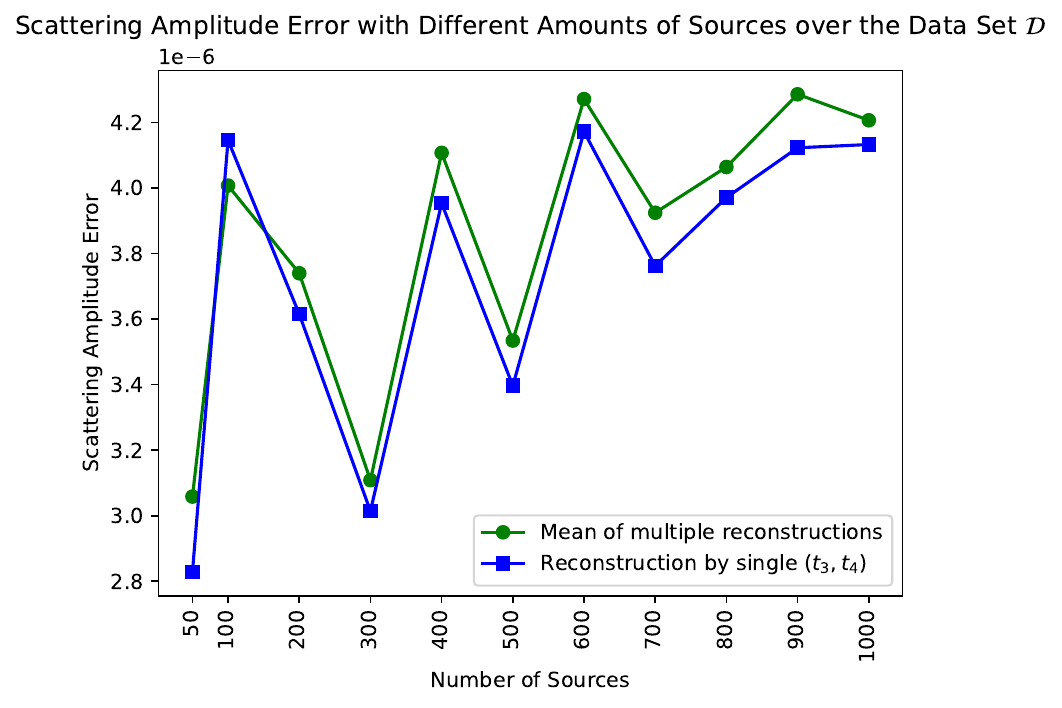}}
		\caption{\footnotesize \it Error behavior with different numbers of sources } \label{fig:x-7-2d}
	\end{figure}

	It is important to determine the optimal number of detectors, sources and $(t_3,t_4)$ for a given experiment. 
	In Figure \ref{fig:x-7-2d} we employ 7 detectors and show the behavior of the error with different numbers of sources on the reconstructions of layer 9. We set $(t_3,t_4)$ equal to the coordinates of 7 detectors, which leads to multiple reconstructions of susceptibility. Then the mean of the multiple reconstructions compared to the single $(t_3,t_4)$ reconstruction is plotted.
	We see that 500 sources almost reaches the minimum of the susceptibility error.  Moreover, the reconstructed susceptibility error by a single $(t_3,t_4)$ is comparable to multiple detectors. Thus, Figure \ref{fig:x-7-2d} shows that one can choose 500 sources and single $(t_3,t_4)$.

\begin{figure}[!htb]
		\centering
		\subfigure{\includegraphics[width=0.3515\textwidth]{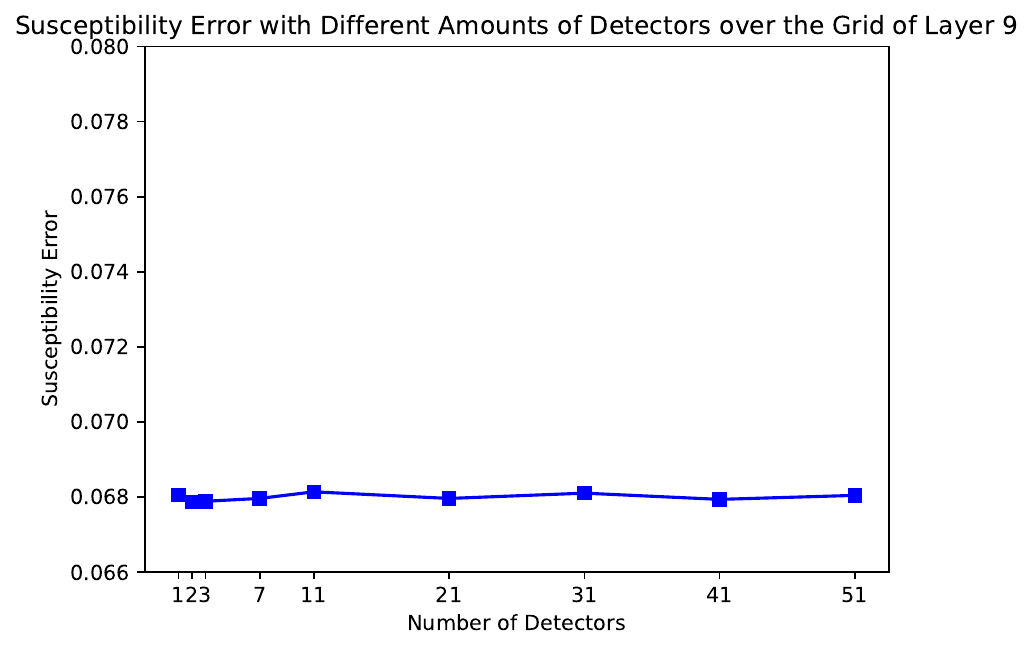}}
		\subfigure{\includegraphics[width=0.3515\textwidth]{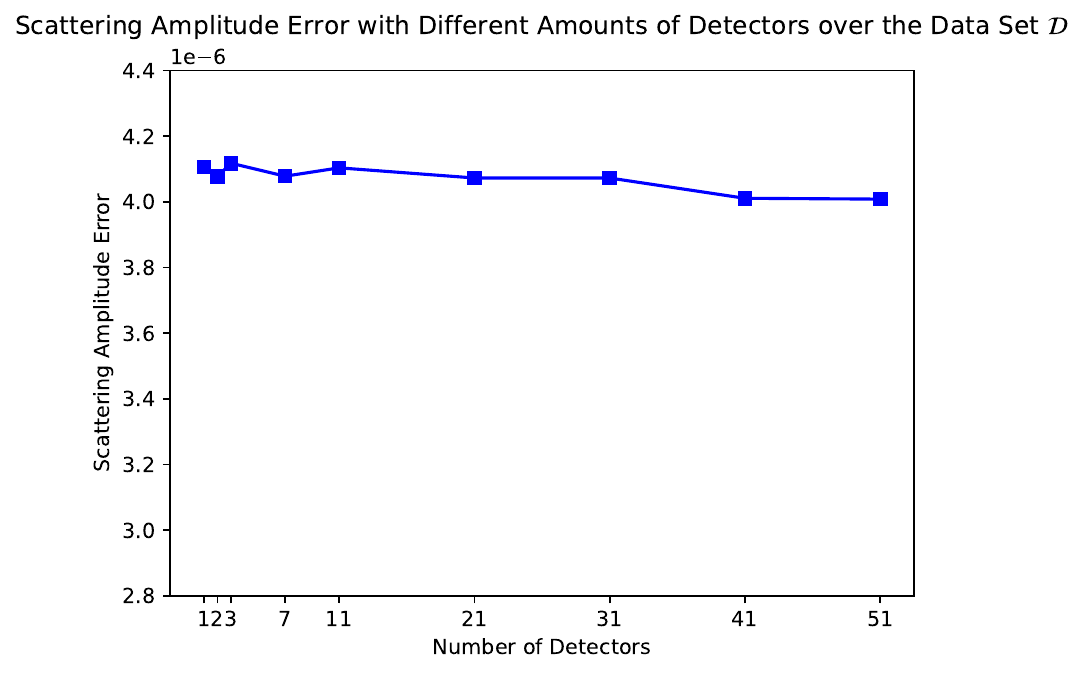}}
		\caption{\footnotesize \it Error behavior using different numbers of detectors} \label{fig:500-x-2d}
	\end{figure}

	Next, in Figure \ref{fig:500-x-2d}, we fix 500 sources and choose a single $(t_3,t_4)$ reconstruction for each grid point, we then test the impact of different numbers of detectors on the reconstructions. We find that the number of detectors does not significantly affect the susceptibility reconstructions. In all  settings (Figure \ref{fig:x-7-2d} and Figure \ref{fig:500-x-2d}),  the scattering amplitude error is on the order of $10^{-6}$. This shows good performance of using the RKHS method.
	Note that the reconstructions depicted in Figure \ref{fig:500-7-2d-slice} have been obtained using the same configuration of 500 sources, 7 detectors, and single $(t_3,t_4)$ as those also utilized in Figure \ref{fig:x-7-2d} and Figure \ref{fig:500-x-2d}.
	
	\subsubsection{Three-dimensional reconstruction method}
\begin{figure}[!htb]
		\centering
		{\includegraphics[width=0.415\textwidth]{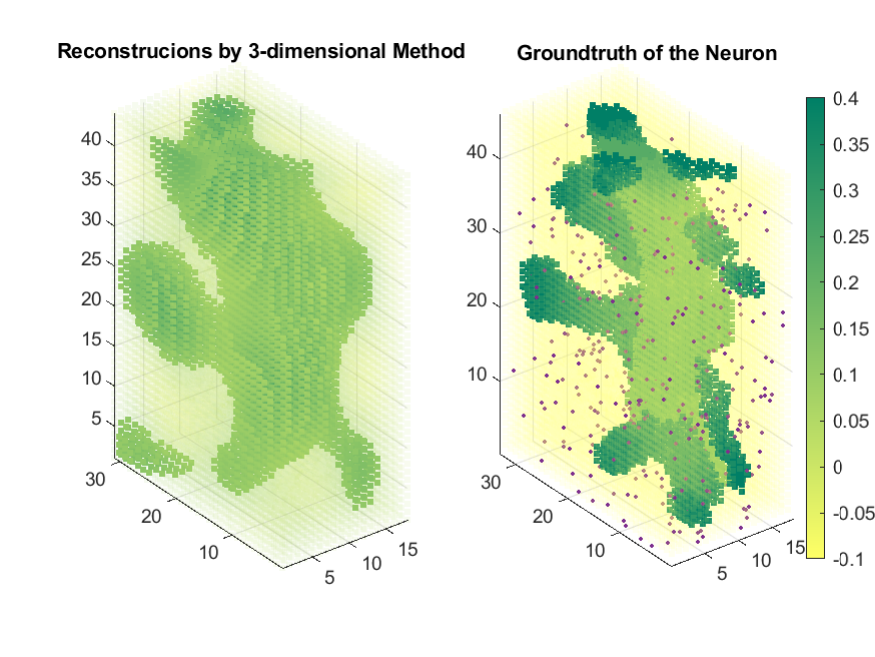}}  
		\caption{\scriptsize \it Neuron reconstructions by the three-dimensional method compared to groundtruth.} \label{fig:500-7-3d-whole}
	\end{figure}

	Here we show the feasibility of our three-dimensional reconstruction method and compare the reconstructions with known results. 
	We segment the neuron into a grid of size 33$\times$46$\times$19, and then we solve the forward problem on this grid with 500 randomly spaced sources and 7 detectors. Subsequently, we reconstruct the model on the grid 
	33$\times$44$\times$18, thus, avoiding an inverse crime.
	The result is shown in Figure \ref{fig:500-7-3d-whole}, where the purple dots indicate the positions of the sources. Once again, we add transparency to the points of value less than 0.08 to better display the colour map.

	We next compared our RKHS method with the method in \cite{GilLevScho18}.
	We consider the same problem setting as in \cite{GilLevScho18} and compute the corresponding scattering amplitude error results in Table \ref{tb:comp-john}.
	We see that for setup in  \cite{GilLevScho18}, the scattering amplitude error by the proposed method is less than reported in \cite{GilLevScho18}.  Moreover, one can reduce the number of detectors from 225 to 7 and still obtain the same scattering amplitude error.
	Thus the proposed RKHS method is superior to the method presented in \cite{GilLevScho18}.
	
	\begin{center}
		\captionof{table}{Comparison between our RKHS method and the one reported in \cite{GilLevScho18}.
		}
		\scalebox{0.7}{
			\begin{tabular}{c|c|c|c}\label{tb:comp-john}
				&  {Linear system \cite{GilLevScho18}} & {RKHS method} & {RKHS method}   \\ 	\hline
				Number of sources & 500 & 500 & 500 \\
				Number of detectors & 225 & 225 & 7 \\
				Discretization of the forward problem & 33$\times$46$\times$19 &33$\times$46$\times$19 &33$\times$46$\times$19   \\
				Scattering amplitude error  & $1*10^{-2}$& $1.9*10^{-4}$ & $1*10^{-4}$ \\			
				\hline
			\end{tabular}
		}
	\end{center}

	
	\section{Conclusion}\label{sec:conclusion}
	
	We have developed a method to reconstruct the dielectric susceptibility of an inhomogeneous scattering medium from internal sources using a RKHS approach. 
	Our numerical results in two- and three dimensions problems exhibit 
	better performance compared to previous work in \cite{GilLevScho18}.  
	We note that the approach we have presented is quite general and applicable to imaging with any scalar wave field with internal sources.

 \subsection*{Acknowledgements}
	This research was funded in whole, or in part, by the Austrian Science Fund (FWF)
	10.55776/P34981 (OS $\&$ YD) – New Inverse Problems of Super-Resolved Microscopy (NIPSUM), SFB 10.55776/F68 (OS) “Tomography Across the Scales,” project F6807-N36 (Tomography with Uncertainties). The financial support by the Austrian Federal Ministry for Digital
	and Economic Affairs, the National Foundation for Research, Technology and Development and the
	Christian Doppler Research Association is gratefully acknowledged.
	The work of JCS was supported by NSF grant DMS-1912821 and AFOSR grant FA9550-19-1-0320.
	For the purpose of open access, the author has applied a CC BY public copyright
	license to any Author Accepted Manuscript version arising from this submission. 
	The computational results presented have been achieved using the Vienna Scientific Cluster (VSC).

	\section*{References}

	\renewcommand{\i}{\ii}
	\printbibliography[heading=none]

	\appendix
	
\end{document}